\input amstex
\input epsf
\documentstyle{amsppt}
\magnification=1200
\NoBlackBoxes
\nologo
\tolerance 1000
\hsize 16.4 truecm
\vsize 22.9 truecm

\global\topskip=0pt
\baselineskip 18pt

\def\rank{\operatorname{rank}}
\def\codim{\operatorname{codim}}
\def\Sym{\operatorname{Sym}}
\def\Sp{\operatorname{Sp}}
\def\Supp{\operatorname{supp}}
\def\gr{\operatorname{gr}}
\def\Arf{\operatorname{Arf}}
\def\var{\operatorname{var}}

\def\V{\Cal V}
\def\K{{\Cal K}}
\def\A{\Cal A}
\def\DD{\Cal D}
\def\W{{\Cal W}}

\def\kk{\frak k}
\def\H{\frak H}
\def\D{\frak D}

\def\G{{\frak G}}

\define \bR {\Bbb R}
\define \bC {\Bbb C}
\define \bZ {\Bbb F}
\define \bF {\Bbb F}

\def\ku{\kappa}
\def\ee{\varepsilon}
\define \la {\lambda}
\define \Ga {\Gamma}
\define \al {\alpha}
\define \siij {\sigma_{i,j}}
\define \Si {\Sigma}
\define \si {\sigma}
\def \de {\delta}
\def \De {\Delta}

\def \ga {h}
\define \pa {\partial}

\def\gs{\geqslant}
\def\ls{\leqslant}

\topmatter
\title Skew-symmetric vanishing lattices and intersections of Schubert 
cells
 \endtitle
\author B.~Shapiro$^{\ddag}$, M.~Shapiro$^*$ and A.~Vainshtein$^{\dag}$
\endauthor
\headline{\hss{\rm\number\pageno}\hss}
\footline{}
\affil $^\ddag$ Department of Mathematics, University of Stockholm\\
S-10691, Sweden,
{\tt shapiro\@matematik.su.se}\\ $^*$ Department of Mathematics,
Royal Institute of Technology\\
S-10044, Sweden, {\tt mshapiro\@math.kth.se} \\ $^\dag$
Dept. of Mathematics and Dept. of Computer Science, University of Haifa\\ 
Mount Carmel, 31905 Haifa, Israel
{\tt alek\@mathcs11.haifa.ac.il} \endaffil




\endtopmatter

\document

\heading {\S 1. Introduction and results} \endheading

In the present paper we apply the theory of skew-symmetric vanishing 
lattices developed around 15 years ago by B.~Wajnryb, S.~Chmutov, and 
W.~Janssen for the necessities of the singularity theory to the 
enumeration of connected components in the intersection of two open 
opposite Schubert cells in the space of complete real flags.
Let us briefly recall the main topological problem considered in \cite {SSV}
and reduced there to a group-theoretical question solved below. Let $N^{n+1}$ 
be the group of real unipotent
uppertriangular $(n+1)\times (n+1)$ matrices and $D_i$ be the
determinant of the submatrix formed by the first $i$ rows and the last 
$i$ columns. Denote by $\Delta_i$ the divisor $\{D_i=0\}\subset N^{n+1}$
and let  $\Delta^{n+1}$ be the union $\cup_{i=1}^n \Delta_i$. Consider now the
complement $U^{n+1}= N^{n+1}\setminus\Delta^{n+1}$. The space $U^{n+1}$ can be
interpreted as the intersection of two open opposite Schubert cells in 
$SL_{n+1}(\bR)/B$.

In \cite {SSV} we have studied the number of connected components in 
$U^{n+1}$. The main result of \cite {SSV} can be stated as follows.

Consider the vector space $T^n(\bZ_2)$ of upper triangular matrices with
$\bF_2$-valued entries. We define the group $\G_n$ as the subgroup of
$GL(T^n(\bF_2))$ generated by $\bZ_2$-linear transformations $g_{ij}$, 
$1\ls i\ls j\ls n-1$. The generator $g_{ij}$ acts on a matrix 
$M\in T^n(\bF_2)$ as follows. Let  $M^{ij}$ denote the $2\times 2$ submatrix
of $M$ formed by  rows $i$ and $i+1$ and columns $j$ and $j+1$
(or its upper triangle in case $i=j$). Then $g_{ij}$
applied to $M$  changes  $M^{ij}$ by adding to each entry of $M^{ij}$ the 
$\bF_2$-valued trace of $M^{ij}$, and does not change all the other entries
of $M$. For example, if $i<j$, then $g_{ij}$ changes $M^{ij}$ as follows:
$$\pmatrix m_{ij}    & m_{i,j+1} \cr
           m_{i+1,j} & m_{i+1,j+1} 
\endpmatrix 
\mapsto 
\pmatrix m_{i+1,j+1} & m_{i,j+1}+m_{ij}+m_{i+1,j+1} \cr
         m_{i+1,j}+m_{ij}+m_{i+1,j+1} & m_{ij}
\endpmatrix.
$$

\noindent All the other entries of $M$ are preserved. The above action
on $T^{n}(\bZ_2)$ is called the {\it first $\G_n$-action}.

\proclaim{Proposition {\rm(Main Theorem of \cite{SSV})}}
The number $\sharp_{n+1}$ of connected components in $U^{n+1}$ coincides
with the number of orbits of the first $\G_n$-action.
\endproclaim

Here we calculate the number of orbits of the first $\G_n$-action
and prove the following result (conjectured in \cite {SSV}).

\proclaim{Main Theorem} 
The number $\sharp_{n+1}$ of connected components in $U^{n+1}$
{\rm(}or{\rm,} equivalently{\rm,} the number of orbits of the first 
$\G_n$-action{\rm)} equals $3\times 2^{n}$ for all $n\gs 5$.
\endproclaim

Cases $n+1=2,3,4 \text{ or }5$ are exceptional and $\sharp_2=2$,
$\sharp_3=6$, $\sharp_4=20$,
$\sharp_5=52$.

The structure of the paper is as follows. In \S 2 we give a detailed
description of the first $\G_n$-action and its quotient by the subspace
of invariants (called the second $\G_n$-action). We formulate explicit
conjectures about the types, number and cardinalities of orbits of both 
actions. In \S 3 we find linear invariants of these actions.
In \S 4 we formulate a number of results about the monodromy group
of a skew-symmetric vanishing lattice over $\bF_2$. Using these results  
we count in  \S 5 the number of orbits of the second $\G_n$-action. Finally,
in \S 6 we explain the relation of the original (first) $\G_n$-action 
to the monodromy group of the corresponding vanishing lattice and 
prove the results stated in \S 2. The concluding \S 7
contains some final remarks and speculations about the origin and further
applications of the above $\G_n$-action in Schubert calculus.

We want to thank Andrei Zelevinsky for a number of stimulating discussions
and his warm support of the whole project. We are indebted to Torsten Ekedahl 
for pointing out the importance of the dual representation of the group 
$\G_n$, to Anatoly Libgober for the key reference \cite{Ja}, and to Yakov
Shlafman for writing sofisticated computer programms that enabled us to collect
enough numerical data for putting forward the main conjecture. 
The third author expresses his gratitude to the Department of 
Mathematics, Royal Institute of Technology in Stockholm, for the 
hospitality and financial support during the final stage of the 
preparation of the manuscript.

\heading \S 2. Two $\G_n$-actions on $\bZ_2$-valued uppertriangular matrices\\
and their orbits. Main tables.
\endheading

\subheading {2.1. The first $\G_n$-action and its invariants}
Let $G$ be a group acting on a linear space $\V$ over $\bF_2$. We say that
$x\in \V$ is
an {\it invariant\/} of the action if $g(x)=x$ for any $g\in G$, and that
$f\in \V^*$ is a {\it dual invariant\/} if  $(g(x),f)=
(x,f)$ for any $x\in \V$, $g\in G$ (here $(\cdot,\cdot)$ is the
standard coupling $\V\times \V^*\to \bF_2$). Evidently, dual 
invariants are just the invariants of the conjugate action of $G$ on $\V^*$.

In what follows we identify $(T^n(\bF_2))^*$ with the space of $\bF_2$-valued 
uppertriangular matrices in such a way that for any pair $M\in T^n(\bF_2)$,
$M'\in (T^n(\bF_2))^*$ one has 
$(M, M')=\sum_{1\ls i\ls j\ls n} m_{ij}m'_{ij}$.

Let us define the matrices $R_i\in (T^n(\bF_2))^*$, $1\ls i\ls n$, as follows:
all the entries of the rectangular
submatrix of $R_i$ formed by the first $i$ rows and the last $n+1-i$ columns
are ones, and all the other entries of $R_i$ are zeros. Thus, for any 
$M\in T^n(\bF_2)$ the value $(M,R_i)$ is just the sum of the entries of
$M \mod 2$ over the corresponding pattern $\rho_i$ (see Fig.~1).

\vskip 20pt
\centerline{\hbox{\hskip 0.2cm \epsfysize=5cm\epsfbox{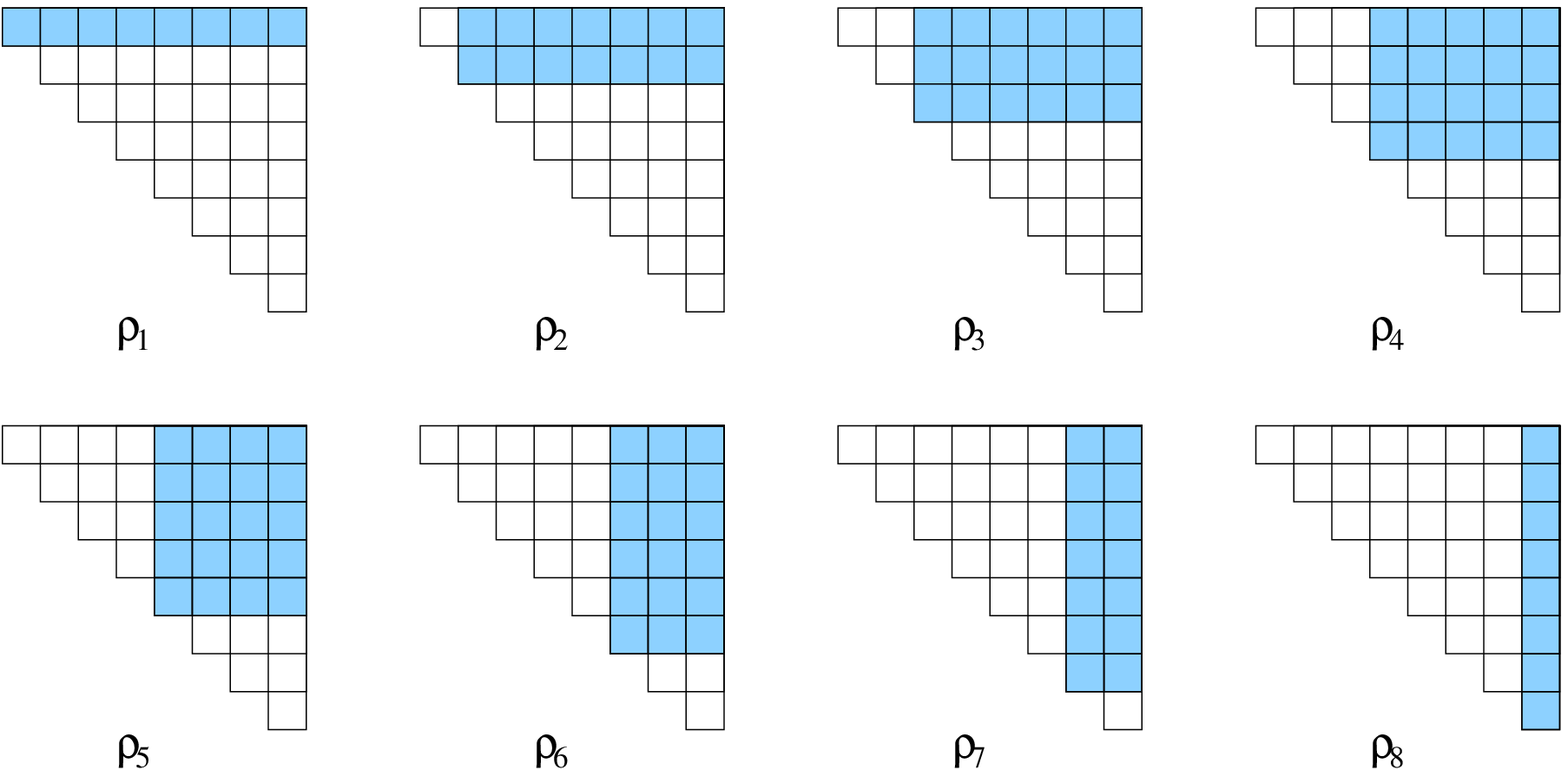}}}
\midspace{1mm} \caption{Fig.~1. Patterns of dual
invariants for the first $\G_n$-action in  case $n=8$}
\vskip 5pt

Let $E_i\in T^n(\bF_2)$, $1\ls i\ls n$, be the matrix whose $i$th diagonal
contains only ones, and all the other entries are zeros.

\proclaim{Theorem 2.1} {\rm (\romannumeral1)} The subspace $\Cal I_n \subset
T^n(\bF_2)$ of invariants of the first $\G_n$-action is an $n$-dimensional
vector space. It has a basis consisting of the matrices $E_1,\dots,E_n$.

{\rm (\romannumeral2)} The subspace $\Cal D_n\subset (T^n(\bF_2))^*$ of dual
invariants of the first $\G_n$-action is an $n$-dimensional vector space. 
It has a basis consisting of the matrices $R_1,\dots,R_n$.
\endproclaim

Let $\Cal D^\perp_n\subset T^n(\bF_2)$ be the subspace  orthogonal to 
$\Cal D_n$ with respect to the standard coupling. The translation of 
$\Cal D^\perp_n$
by an arbitrary element $M\in T^n(\bF_2)$ we call a {\it stratum\/} of
the first $\G_n$-action. By the definition, all the dual invariants are
fixed at all the elements of a stratum. An $n$-dimensional vector 
$h^S=(h_1^S,\dots,h_n^S)$ is said to be the {\it height\/} of a stratum $S$ 
(with respect to the basis $\{R_i\}$)
if $(M,R_i)=h^S_i$ for any $M\in S$. A stratum is called {\it symmetric\/}
if its height is symmetric with respect to its middle, that is, if 
$h_i^S=h_{n-i+1}^S$ for all $1\ls i\ls n$.

Evidently, each stratum is a union of certain orbits of the first 
$\G_n$-action. The structure of strata is described by Theorem~2.2 below. 
For the sake of simplicity, we omit the words ``of the first $\G_n$-action''
in the formulation and write just ``orbit'' and ``stratum''.

\proclaim{Theorem 2.2} {\rm (\romannumeral1)} Let $n=2k+1\gs5${\rm,} then

each of $2^{k+1}$ symmetric strata consists of one orbit of length
$2^{2k^2+k-1}-\ee_k2^{k^2+k-1}${\rm,} one orbit of length 
$2^{2k^2+k-1}+\ee_k2^{k^2+k-1}-2^k${\rm,} and $2^k$ orbits of length $1${\rm,}
where $\ee_k=-1$ for $k=4t+1$ and $\ee_k=1$ otherwise{\rm;}

each of $2^n-2^{k+1}$ nonsymmetric strata consists of two orbits
of length $2^{2k^2+k-1}$.

{\rm (\romannumeral2)} Let $n=2k\gs6${\rm,} and let $\bar h$ denote the 
vector of 
length $n$ the first $k$ entries of which are equal to $(1,0,1,0,...)$ and 
the last $k$ entries vanish{\rm,} then

each of $2^{k}$ symmetric strata consists of two orbits of length
$(2^{2k(k-1)}-1)2^{k-1}$ and $2^k$ orbits of length $1${\rm;} 

each of $2^k$ nonsymmetric strata $S$ such that $h^S-\bar h$
is symmetric with respect to its middle consists of one orbit of length
$(2^{k(k-1)}-1)2^{k^2-1}$ and one orbit of length 
$(2^{k(k-1)}+1)2^{k^2-1}${\rm;}

each of $2^n-2^{k+1}$ remaining nonsymmetric strata consists of two 
orbits of length $2^{2k^2-k-1}$.
\endproclaim

The assertions of Theorem~2.2 can be summarized in the following table.
Observe that the total number of orbits in both cases equals $3\times2^n$,
and we thus get Main Theorem.

$$\matrix
&	\vdots	&n=2k+1	&	& \vdots
&	n=2k	&	& \vdots\\
\hdotsfor{8}\\
\text{type}& \vdots & \text{cardinality}
& \sharp_{orb} &\vdots & \text{cardinality} & \sharp_{orb} &\vdots\\
\hdotsfor{8}\\
\text {trivial}& \vdots &	1	& 2^{2k+1}&\vdots&1 &
2^{2k} &\vdots \\
\text{standard}&\vdots &2^{2k^2+k-1} & 2^{2k+2}-2^{k+2}&\vdots & 2^{2k^2-k-1}
&2^{2k+1}-2^{k+2}&\vdots \\ 
\text {type 1}&\vdots &(2^{k^2}-\ee_k) 2^{k^2+k-1}
&2^{k+1}&\vdots & 0&0&\vdots \\
\text {type 2}&\vdots &(2^{k^2}-\ee_k) (2^{k^2-1}+\ee_k) 2^{k}&2^{k+1}&\vdots &
0&0&\vdots \\
\text {type 3}&\vdots &0&0&\vdots &
(2^{2k(k-1)}-1) 2^{k-1}&2^{k+1}&\vdots \\ \text {type 4}&\vdots &0&0&\vdots &
(2^{k(k-1)}-1)2^{k^2-1}&2^k&\vdots \\
\text {type 5}&\vdots &0&0&\vdots &
(2^{k(k-1)}+ 1) 2^{k^2-1}&2^k&\vdots \\
\hdotsfor{8}
\endmatrix $$

\centerline {Table 1. The orbits of the first $\G_n$-action} 
\bigskip

\subheading{2.2. The second $\G_n$-action and its invariants}
Let us introduce a $\G_n$-action on $T^{n-1}(\bZ_2)$ closely related
to the first $\G_n$-action on $T^n(\bZ_2)$, which we will call the 
{\it second $\G_n$-action}. This action is induced by taking the quotient 
modulo the subspace $\Cal I_n$ of the invariants of the first $\G_n$-action.
Recall that by Theorem~2.1(i) $\Cal I_n$ is an $n$-dimensional
subspace of all matrices having the same entry on each diagonal. One can 
suggest the following natural description of the second $\G_n$-action.

Consider the linear map $\Psi_n\: T^n(\bZ_2)\to T^{n-1}(\bZ_2)$
such that the $(i,j)$th entry in the image equals the sum of the $(i,j)$th
and the $(i+1,j+1)$th entries in the inverse image (thus entry $(i,j)$ in the
image can be considered as representing the submatrix $M^{ij}$ of the initial
matrix). Evidently, $\ker \Psi_n=\Cal I_n$, and we obtain the induced
action of $\G_n$ on $T^{n-1}(\bZ_2)$. For any $M\in T^{n-1}(\bZ_2)$ the 
generator $g_{ij}$ affects only the $3\times 3$ submatrix of $M$ centered  
at $m_{ij}$ by adding $m_{ij}$ to all entries marked by asterisks in the 
shape below:
$$
\pmatrix *&*&{}\\*&m_{ij}&*\\{}&*&*\endpmatrix.\tag1
$$
Here we use the convention that if the above $3\times 3$-shape does not
fit completely in the upper triangle, then we change only the 
entries that fit.

For any $i$, $1\ls i\ls k=[n/2]$, we define a subset $\pi_i$ of the entries 
of the $(n-1)\times(n-1)$ triangular shape as follows.
The subset $\pi_k$ coincides with the initial $(n-1)\times(n-1)$ triangular 
shape. To build $\pi_i$, $1\ls i\ls k-1$, cut off the three
$i\times i$ triangular shapes placed in all the three corners of the initial
shape and get a hexagonal shape $\chi_i^1$ with the sides of lengths $i+1$
and $n-2i-1$. The shape $\chi_i^1$ is included into a nested family of
$j_i=\min\{i+1,n-2i-1\}$ hexagonal shapes $\chi_i^j$; the shape $\chi_i^j$ is
obtained by peeling the external layer of $\chi_i^{j-1}$. The last shape 
$\chi_i^{j_i}$ degenerates to a triangle (upper if $j_i=i+1$ and lower 
otherwise). The subset $\pi_i$ consists of the three $i\times i$ triangular 
shapes as above and $[j_i/2]$ layers of the form $\chi_i^{2j}\setminus 
\chi_i^{2j+1}$, $1\ls j\ls[j_i/2]$; here we assume that 
$\chi_i^{j_i+1}=\varnothing$, and so for $j_i$ even the last layer 
consists of the whole triangle $\chi_i^{j_i}$. 

We now define the matrices $P_i\in (T^{n-1}(\bF_2))^*$, $1\ls i\ls k$, as 
follows: all the entries of $P_i$ that belong to $\pi_i$ are ones, and all the
other entries of $P_i$ are zeros. Thus, for any 
$M\in T^{n-1}(\bF_2)$ the value $(M,P_i)$ is just the sum of the entries of
$M \mod 2$ over the corresponding pattern $\pi_i$ (see Fig.~2).

\vskip 20pt
\centerline{\hbox{\hskip 0.2cm \epsfysize=3cm\epsfbox{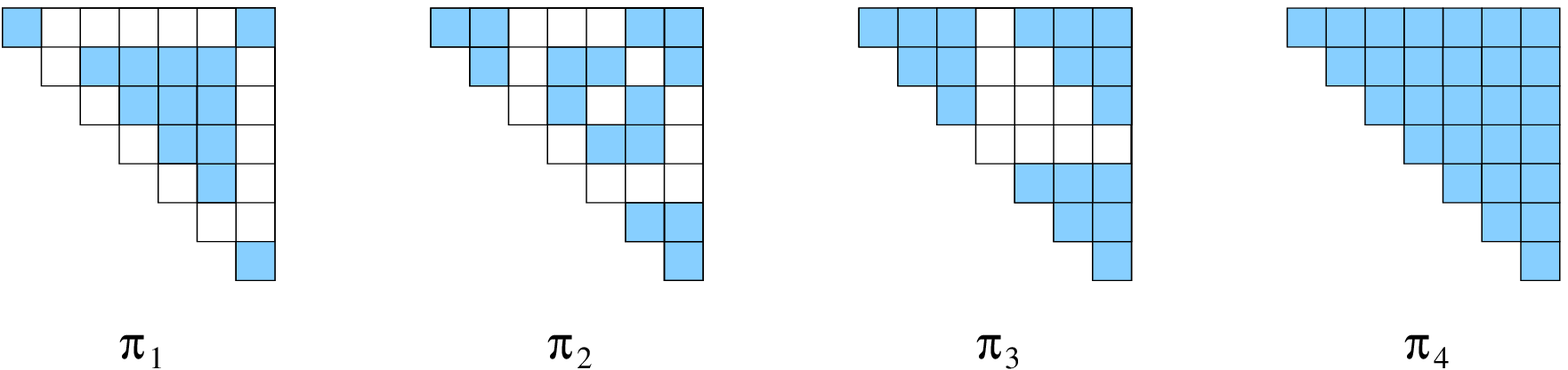}}}
\midspace{1mm} \caption{Fig.2. Patterns of dual invariants for the 
second $\G_n$-action in case $n=8$}
\vskip 5pt

\proclaim{Theorem 2.3} {\rm (\romannumeral1)} The subspace 
of invariants of the second $\G_n$-action is trivial.

{\rm (\romannumeral2)} The subspace $\D_{n-1}\subset (T^{n-1}(\bF_2))^*$ 
of dual invariants of the second $\G_n$-action is a $k$-dimensional vector 
space{\rm,} $k=[n/2]$. It has a basis consisting of the matrices 
$P_1,...,P_k$.
\endproclaim

It follows immediately from Theorem~2.3(ii) that $\widetilde P_i=P_i+P_{i-1}$,
where $P_0=0$, is also a basis of  $\D_{n-1}$.  
The strata of the second $\G_n$-action and their heights (with respect
to the basis $\{\widetilde P_i\}$), which we denote $\eta$, are defined 
in the same way as for the first $\G_n$-action. It is easy to see that 
$\Psi_n$ takes any stratum 
of the first $\G_n$-action to a stratum of the second $\G_n$-action.
Denote by $\psi_n\:\bF_2^n\to\bF_2^k$ the mapping that takes $h^S$ to
$\eta^{\Psi_n(S)}$. The following statement gives an explicit description
of $\psi_n$.  

\proclaim{Theorem 2.4} The heights of the strata $S$ and $\Psi_n(S)$  
with respect to the bases $\{R_i\}$ and $\{\widetilde P_i\}$
satisfy the relations
$$\gather
\eta^{\Psi_n(S)}_i=h^S_i+h^S_{i+1}+h^S_{n-i}+h^S_{n-i+1}
\quad\text{for}\quad 1\ls i<k,\\
\eta^{\Psi_n(S)}_k=h^S_k+h^S_{n-k+1}.
\endgather
$$
\endproclaim

The structure of the strata of the second $\G_n$-action
is described in Theorem~2.5 below.
For the sake of simplicity, we omit the words ``of the second $\G_n$-action''
in the formulation and write just ``orbit'' and ``stratum''.

\proclaim{Theorem 2.5} {\rm (\romannumeral1)} Let $n=2k+1\gs5${\rm,} then

the stratum at height $(0,\dots,0)$ consists of
one orbit of length $2^{2k^2-1}-\ee_k2^{k^2-1}${\rm,}
one orbit of length $2^{2k^2-1}+\ee_k2^{k^2-1}-1${\rm,}
and one orbit of length $1${\rm,} where $\ee_k=-1$ for $k=4t+1$ and 
$\ee_k=1$ otherwise{\rm;} 

each of the other $2^k-1$ strata is an orbit of length $2^{2k^2}$.

{\rm (\romannumeral2)} Let $n=2k\gs6${\rm,} and let $\bar\eta$ denote the 
vector of length $k$ the first $k-1$ entries of which are equal to $1$ and 
the last entry equals $k \mod 2${\rm,} then

the stratum at height $(0,\dots,0)$ consists of one orbit of length
$2^{2k(k-1)}-1$ 
and one orbit of length $1${\rm;}

the stratum at height $\bar\eta$ consists of one orbit of length 
$(2^{k(k-1)}-1)2^{k(k-1)-1}$ 
and one orbit of length $(2^{k(k-1)}+1)2^{k(k-1)-1}${\rm;} 

each of the other $2^k-2$ strata is an orbit of length $2^{2k(k-1)}$. 
\endproclaim

The assertions of Theorem~2.5 can be summarized in the following table.
Observe that the total number of orbits in both cases equals $2^k+2$.

$$\matrix
&	\vdots	&
n=2k+1 & 	& \vdots
&	n=2k	&		& \vdots\\
\hdotsfor{8}\\
\text{type}& \vdots & \text{cardinality}
& \sharp_{orb} &\vdots & \text{cardinality} & \sharp_{orb}&\vdots\\
\hdotsfor{8}\\
\text {trivial}& \vdots &	1	&1&\vdots&1 & 1&\vdots \\
\text {standard}&\vdots &2^{2k^2}&2^k-1&\vdots &
2^{2k(k-1)}&2^{k}-2&\vdots \\
\text {type 1}&\vdots &(2^{k^2}-\ee_k)2^{k^2-1}&1&\vdots & 0&0&\vdots \\
\text {type 2}&\vdots &(2^{k^2}-\ee_k)(2^{k^2-1}+\ee_k)&1&\vdots&0&0&\vdots \\
\text {type 3}&\vdots &0&0&\vdots &(2^{2k(k-1)}-1)&1&\vdots \\
\text {type 4}&\vdots &0&0&\vdots &(2^{k(k-1)}-1)2^{k(k-1)-1}&1&\vdots \\ 
\text{type 5}&\vdots &0 &0&\vdots &(2^{k(k-1)}+1)2^{k(k-1)-1} &1&\vdots \\
\hdotsfor{8}
\endmatrix $$

\centerline {Table 2. The orbits of the second $\G_n$-action}
\bigskip

\remark{Remark}  The map $\Psi_n$ sends the orbits of the first $\G_n$-action 
to the orbits of the second $\G_n$-action with the same name (i.e. trivial 
to trivial, standard to standard, etc). \endremark

\heading \S 3. Invariants and dual invariants of the $\G_n$-actions 
\endheading

In this section we prove Theorems~2.1,~2.3 and~2.4.

\subheading {3.1. The first $\G_n$-action}
Recall that we have identified $(T^n(\bF_2))^*$ with the space of 
$\bF_2$-valued 
uppertriangular matrices in such a way that for any pair $M\in T^n(\bF_2)$,
$M'\in (T^n(\bF_2))^*$ one has 
$(M, M')=\sum_{1\ls i\ls j\ls n} m_{ij}m'_{ij}$.

\proclaim{Lemma 3.1} The conjugate to the first $\G_n$-action is given by 
$$
(g_{ij}M')_{kl}=
\left\{
\alignedat2
&m'_{i+1,j}+m'_{i,j+1}+m'_{i+1,j+1}\qquad &&\text{if $k=i$, $l=j$}, \\
&m'_{i+1,j}+m'_{i,j+1}+m'_{i,j}\qquad &&\text{if $k=i+1$, $l=j+1$},\\
&m'_{k,l}\qquad &&\text{otherwise.}
\endalignedat
\right.
$$
\endproclaim

\demo{Proof} Indeed, the action of $g_{ij}$ on $T^n(\bF_2))$ affects only 
$M^{ij}$, while the action of $g_{ij}$ on $(T^n(\bF_2))^*$ defined above
affects only $(M')^{ij}$. Therefore, it suffices to consider the actions in the
corresponding 4-dimensional spaces. The matrix of $g_{ij}$ in coordinates
$m_{ij}$, $m_{i,j+1}$, $m_{i+1,j}$, $m_{i+1,j+1}$ is
$$
\pmatrix 0 & 0 & 0 & 1 \cr
         1 & 1 & 0 & 1 \cr
         1 & 0 & 1 & 1 \cr
         1 & 0 & 0 & 0 \endpmatrix,
$$
hence its conjugate is exactly as asserted by the lemma.
\qed
\enddemo

Now we are ready to prove Theorem~2.1.

\demo{Proof of Theorem~2.1} (i) Follows immediately from the fact that
$\Cal I_n$ is defined by the equations $m_{ij}+m_{i+1,j+1}=0$,
$1\ls i\ls j\ls n-1$.

(ii) By Lemma~3.1, $\Cal D_n$ is defined by the equations $m'_{ij}+
m'_{i,j+1}+m'_{i+1,j}+m'_{i+1,j+1}=0$, $1\ls i\ls j\ls n-1$ (with the same
as above convention concerning the case when only a part of the submatrix
fits into the uppertriangular shape). One can prove easily that these
$\binom {n-1}{2}$ equations are linearly independent. Indeed, each of
$m'_{1j}$, $1<j<n$, enters exactly two equations, while $m'_{11}$ only
one equation. It follows immediately that none of the equations involving
$m'_{1j}$ may participate in a nontrivial linear combination. The rest
of the equations correspond to the same situation in dimension $n-1$, so
their linear independence follows by induction. We thus get $\dim\Cal D_n=n$.
It remains to show that $R_1,\dots,R_n$ provide a basis for $\Cal D_n$.
Evidently, all these matrices are linearly independent and
satisfy the equations defining $\Cal D_n$.
\qed
\enddemo

\subheading {3.2. The second $\G_{n}$-action}
For any matrix entry $m_{ij}$ we define its {\it neighbors\/} as all the
entries marked by asterisks in (1), or, more precisely, those of
$m_{i-1,j-1}$, $m_{i-1,j}$, $m_{i,j-1}$, $m_{i,j+1}$, $m_{i+1,j}$, and 
$m_{i+1,j+1}$ that fit in the uppertriangular shape. It is helpful to
consider a graph $\H_{n-1}$ whose vertices are all the matrix entries and 
edges join each entry with its neighbors. It is easy to see that $\H_{n-1}$
is an equilateral triangle with the sides of length $n-2$ on the triangular 
lattice (see Fig.~3). One can give the following 
description of the second $\G_n$-action in terms of $\H_{n-1}$: $g_{ij}$
acts on the space of $\bF_2$-valued functions on $\H_{n-1}$ by adding the
value at $m_{ij}$ to the values at all of its neighbors. 

\vskip 20pt
\centerline{\hbox{\hskip 0.2cm \epsfysize=3cm\epsfbox{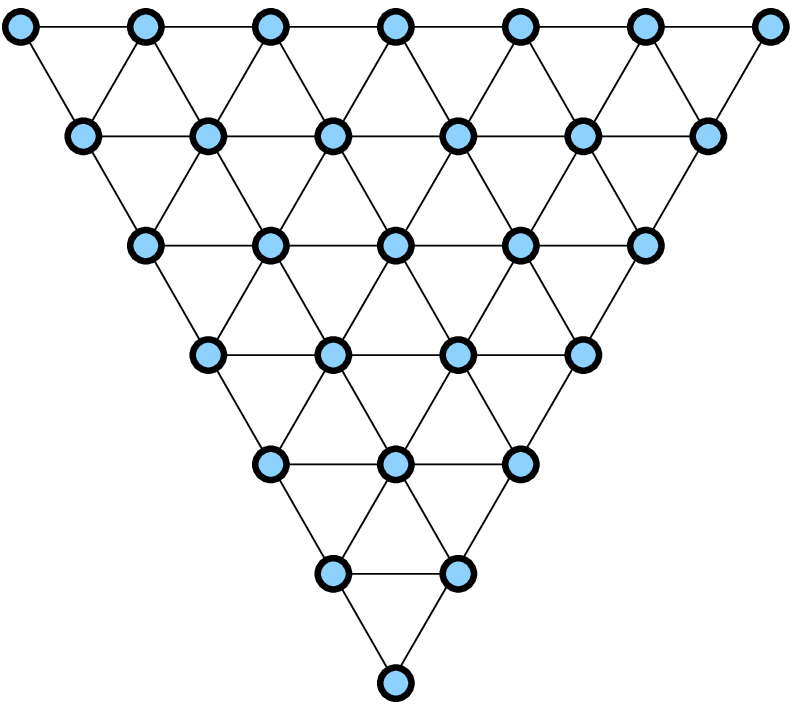}}}
\midspace{1mm} \caption{Fig.3. The graph $\H_{n-1}$ in case $n=8$}
\vskip 5pt

\proclaim{Lemma 3.2} The conjugate to the second $\G_n$-action acts as 
follows{\rm:} $g_{ij}$ adds to the value at $m_{ij}$ the sum of the values at 
all its neighbors.
\endproclaim

\demo{Proof} Similarly to the proof of Lemma~3.1, it suffices to consider
the action in the corresponding 7-dimensional subspaces. The details are
straightforward.
\qed\enddemo

To prove Theorem~2.3, we need  the following technical proposition. 
Consider the subspace $\D^1_{n-1}\subset (T^{n-1}(\bF_2))^*$ of all
linear forms invariant under the subgroup of $\G_n$ generated by 
$\{g_{ij}\: i\gs2\}$; evidently, $\D_{n-1}\subseteq \D^1_{n-1}$.
Let $\omega\:\D^1_{n-1}\to\bF_2^{n-1}$ denote the projection on the 
first row, and let $\Sym^{p}\subset \bF_2^p$ denote the space of all
vectors symmetric with respect to their middle.

\proclaim{Lemma~3.3} {\rm (\romannumeral1)} The $i$th row of an arbitrary
matrix $M\in\D^1_{n-1}$ belongs to $\Sym^{n-i}${\rm,} $1\ls i\ls n-1$.

{\rm (\romannumeral2)} $\dim\D^1_{n-1}=n-1$.

{\rm (\romannumeral3)} The image of $\omega$ coincides with $\Sym^{n-1}$.
\endproclaim

\demo{Proof} For $n=3$ an immediate check shows that $\D^1_2$ consists of
the following four matrices:
$$
\pmatrix 
    0&0\\ 
\quad&0\endpmatrix,\quad
\pmatrix 
    0&0\\ 
\quad&1\endpmatrix,\quad 
\pmatrix 
    1&1\\
\quad&0\endpmatrix,\quad
\pmatrix 
    1&1\\
\quad&1\endpmatrix,
$$
and so all the assertions of the lemma hold true.

Consider now an arbitrary $M\in\D^1_{n-1}$. Let
$a=(a_1,\dots,a_{n-1})$, $b=(b_1,\dots,b_{n-2})$, $c=(c_1,\dots,c_{n-3})$
be the first three rows of $M$. It follows easily from the
definition of $\D^1_{n-1}$ that $a$, $b$, $c$ satisfy equations
$a_i=a_1+b_1+b_i+b_{i-1}+c_{i-1}$, $1\ls i\ls n-1$, where $b_0=c_0=0$.
Evidently, the submatrix of $M$ obtained by deleting the first row
belongs to $\D^1_{n-2}$. By induction, we may assume that the rows of this
matrix are symmetric, that is, $b_i=b_{n-1-i}$, $1\ls i\ls n-2$, and
$c_i=c_{n-2-i}$, $1\ls i\ls n-3$. We thus have $a_{n-i}=a_1+b_1+b_{n-i}+
b_{n-1-i}+c_{n-1-i}=a_i$, $1\ls i\ls n-1$, and hence the first row of
$M$ is symmetric as well. Moreover, it is clear that the solutions of 
the above equations form a 1-dimensional affine subspace, and hence
$\dim\D^1_{n-1}=\dim\D^1_{n-2}+1=n-1$. To get the third statement it 
suffices to notice that the first two rows define $M$ uniquely, and that 
the dimension of the subspace spanned by these two rows is exactly $n-1$.
\qed
\enddemo

\demo{Proof of Theorem~2.3} (i) Follows immediately from the definition.

(ii) Let us define the linear mapping $\nu\:\D^1_{n-1}\to\bF_2^{n-1}$
by the following rule: $\nu_i(M)$ is the sum of the neighbors of $m_{1i}$.
By definition, $\D_{n-1}=\ker\nu$, so in order to find $\dim\D_{n-1}$ it
suffices to determine the image of $\nu$. As in the proof of Lemma~3.3,
denote by $a$ and $b$ the first two rows of $M$. Clearly,
$\nu_i=a_{i+1}+a_{i-1}+b_i+b_{i-1}$, $1\ls i\ls n-1$, where $a_0=a_n=b_0=0$.
These equations define a linear mapping that takes the space $\Sym^{n-1}$ 
of the first rows to the image of $\nu$. It is easy to see that the corank
of this mapping is zero for $n-1=2k$ and one for $n-1=2k-1$; in other words,
it equals $[n/2]-[(n-1)/2]$. Since $\dim\Sym^{n-1}=[n/2]$, we get that the 
dimension of the image of $\nu$ equals $[(n-1)/2]$. It follows immediately
that $\dim\D_{n-1}=n-1-[(n-1)/2]=[n/2]$.

The remaining assertion about the structure of dual invariants $P_{i}$ is 
rather obvious. Its proof consists of the two simple checks: that
$P_{i}$'s are linearly independent and that any matrix entry has an even 
number of nonzero neighbors (shown by shaded areas on Fig.~2).
Both are straightforward.
\qed
\enddemo

Now we can prove Theorem~2.4.

\demo{Proof of Theorem~2.4} It follows easily from Theorem~2.1 that the 
space $D^n(\bF_2)$ of diagonal $n\times n$ matrices is transversal to the 
strata of the first action; moreover, $\Psi^n$ takes $D^n(\bF_2)$ to
$D^{n-1}(\bF_2)$. Therefore it is enough to check the assertion of the theorem
only for diagonal matrices. For $M\in D^n(\bF_2)$ and $M'=\Psi_n(M)$
Theorems~2.1 and~2.3 give
$$\multline
(M',\widetilde P_i)=m'_{ii}+m'_{n-i,n-i}+m'_{n-i+1,n-i+1}=\\
(m_{ii}+m_{i+1,i+1})+(m_{n-i,n-i}+m_{n-i+1,n-i+1})=\\
(M,R_i)+(M,R_{i+1})+(M,R_{n-i})+(M,R_{n-i+1})
\endmultline
$$
for $1\ls i<k$. For $n=2k+1$ this equality holds also for $i=k$,
but now $k+1=n-k$, hence the two middle terms vanish and thus 
$(M',\widetilde P_k)=(M,R_k)+(M,R_{n-k+1})$. Finally, for $n=2k$ we have 
$(M',\widetilde P_k)=m'_{kk}=(M,R_k)+(M,R_{k+1})=(M,R_k)+(M,R_{n-k+1})$.
\qed
\enddemo

\heading \S 4. Some results about skew-symmetric vanishing lattices \endheading

In this section we quote and prove a number of results concerning the
natural action of a
group generated by transvections preserving a given skew-symmetric bilinear
form on a vector space over $\bF_2$. Our main reference is \cite {Ja}.
The relation to our original problem is explained
in details in the next section and is based on the fact that the 
conjugate to the second $\G_n$-action is exactly of this kind.
This allows us to describe completely its orbits, as well as the orbits of
the second $\G_n$-action itself.

\subheading{4.1. Vanishing lattices and their monodromy groups}
We assume  that $\V$ is a vector space over $\bF_2$ equipped with a
skew-symmetric bilinear form $\langle\cdot,\cdot\rangle=(\cdot,
L(\cdot))$, where $L$ is a linear map $L\:\V\to \V^*$.
A {\it  quadratic function\/} $q$ associated with  $\langle
\cdot, \cdot\rangle$ is an arbitrary $\bF_2$-valued function on $\V$
satisfying
$$
q(\lambda x+\mu y)=\lambda^2 q(x)+\mu^2 q(y)+\lambda\mu\langle x,y\rangle.
$$
It is clear that a quadratic function completely determines the
corresponding bilinear form. With any basis $B$ of $\V$ we associate a unique 
quadratic function $q_B$ by requiring it to take value 1 on all elements of 
$B$.

Let $\K=\ker L$ be the kernel of $\langle\cdot,\cdot\rangle$, $\ku=\dim\K$, and
$(e_1,f_1,...,e_m,f_m,g_1,...,g_\ku)$ be a symplectic basis for 
$\langle \cdot, \cdot \rangle$, that is,
$\langle x,y\rangle=\sum (x_iy_i^\prime+y_ix_i^\prime)$, where $x=\sum x_i
e_i+\sum x_i^\prime f_i +
\sum x_i^{\prime \prime}g_i$ and $y=\sum y_i e_i+\sum y_i^\prime f_i +
\sum y_i^{\prime\prime}g_i$.  If the restriction $q|_\K$ vanishes, then one
can define the {\it Arf invariant\/} of $q$ by
$\Arf(q)=\sum_{i=1}^mq(e_i)q(f_i)$, see e.g. \cite{Pf}. For fixed
dimensions of $\V$ and $\K$ there exist, up to isomorphisms,
at most three possibilities: (i) $q(\K)=0$, $\Arf(q)=1$;
(ii) $q(\K)=0$, $\Arf(q)=0$; (iii) $q(\K)=\bF_2$ (and so $\ku\gs 1$).

In the first case one has $\vert q^{-1}(1)\vert =2^{2m+\ku-1}+2^{m+\ku-1}$,
$\vert q^{-1}(0)\vert =2^{2m+\ku-1}-2^{m+\ku-1}$, in the second case one has
$\vert q^{-1}(1)\vert =2^{2m+\ku-1}-2^{m+\ku-1}$,
$\vert q^{-1}(0)\vert =2^{2m+\ku-1}+2^{m+\ku-1}$, and in the third case one has
$\vert q^{-1}(1)\vert=\vert q^{-1}(0)\vert =2^{2m+\ku-1}$.

For any $\de\in \V$ we define the  {\it symplectic transvection\/}
$T_\de\:\V\to \V$ by $T_\de(x)=x- \langle x, \de \rangle \de$. (Notice
that $T_\de$ is an element of the group $\Sp^\sharp \V$ of the 
automorphisms of $(\V,\langle\cdot,\cdot \rangle)$).
Given any subset $\Delta\in V$, we let $\Ga_\Delta\subseteq \Sp^\sharp \V$
denote the subgroup generated by the transvections $T_\delta$, $\delta\in
\Delta$.

The main object of this section is a  {\it vanishing lattice}, that is,
a triple $(\V,\langle\cdot,\cdot\rangle,\Delta)$ satisfying the
following three conditions:
(i) $\Delta$ is a $\Ga_\Delta$-orbit; (ii) $\Delta$ generates $\V$;
(iii) if $\rank \V>1$, then there exist $\delta_1,\delta_2\in \Delta$ such
that $\langle \delta_1,\delta_2 \rangle=1$. 
The group $\Ga_\De$ is called the {\it monodromy group\/}
of the lattice.

We say that a basis $B$ of $\V$
is {\it weakly distinguished\/} if $\Ga_B=\Ga_\Delta$. In this case
$\Ga_B$ respects $q_B$, so, in particular, $q_B(\delta)=1$ for all
$\delta \in \Delta$.
Bases $B$ and $B^\prime$ are called {\it equivalent\/} if $B^\prime\subset
\Ga_B\cdot B$ and
$B\subset \Ga_{B^\prime}\cdot B^\prime$. One can easily see that if $B$ and
$B^\prime$
are equivalent, then $\Ga_B=\Ga_{B^\prime}$, $q_B=q_{B^\prime}$, and if $B$
is a weakly
distinguished basis for $(\V,\langle\cdot,\cdot\rangle,\Delta)$, then so is
$B^\prime$.

Let $B=(b_1,...,b_d)$ be a basis in $\V$, $d=\dim \V$. 
We define the graph $\gr(B)$
of $B$ as follows: $B$ is its vertex set, and $b_i$ is connected
by an edge with $b_j$ if $\langle b_i,b_j \rangle =1$. 
It is easy to see that $\gr(B)$ is connected if $B$
is weakly distinguished.

A basis is called {\it special\/} if it is equivalent to a basis
$B=(b_1,...,b_d)$
such that for some $k$, $1\ls k\ls d$, we have $\langle b_i,b_j\rangle=0$ iff
$i=j$ or $i,j\gs k+1$, and {\it nonspecial\/} otherwise.
A vanishing lattice admitting a special (resp. nonspecial) weakly distinguished
basis and its monodromy group are
called {\it special\/} (resp. {\it nonspecial\/}).  

Nonspecial monodromy groups have an especially simple characterization
(we are primarily interested in this case since for $n\gs 5$ the
group $\G_n$ introduced in \S 1 can be interpreted as a nonspecial monodromy 
group, see \S 5 for details.)

\proclaim{Theorem~4.1 {\rm (\cite {Ja}, Th.~3.8)}} 
Let $(\V,\langle \cdot,\cdot\rangle,\Delta)$ be a
vanishing lattice admitting a nonspecial weakly distinguished basis $B$.
Then $\Ga_B$
coincides with the subgroup $O^\sharp(q_B)$ of $\Sp^\sharp V$ consisting of
all automorphisms that preserve $q_B$.
\endproclaim

It follows that the classification of nonspecial vanishing lattices reduces
to the classification of the corresponding quadratic functions. We thus see
that for given $m$ and $\ku$ such that $\dim \V=2m+\ku$, $\dim \K=\ku$ there
exist exactly three nonspecial vanishing lattices, depending on the values
of $q_B(\K)$ and $\Arf(q_B)$. They are denoted by $O_1^\sharp(2m,\ku,\bF_2)$,
$O_0^\sharp(2m,\ku,\bF_2)$, and $O^\sharp(2m,\ku,\bF_2)$, and correspond to
the cases $\Arf(q_B)=1$, $\Arf(q_B)=0$, and $q_B(\K)=\bF_2$, respectively
(see \cite{Ja, 4.2} for details).

The following statement, which can be extracted easily from \cite{Ja, \S4},
provides a sufficient condition for a vanishing lattice to be nonspecial.

\proclaim{Lemma~4.2} A vanishing lattice is nonspecial if
it admits a weakly distinguished basis $B$ such that $\gr(B)$ contains
the standard Dynkin diagram of the Coxeter group $E_6$ as an induced subgraph.
\endproclaim

\subheading{4.2. Orbits of nonspecial monodromy groups and of the conjugate 
actions}
First of all, let us find the number of orbits of the monodromy group in
the nonspecial case.

\proclaim{Lemma~4.3} Let $(\V,\langle \cdot,\cdot\rangle,\Delta)$ be a
vanishing lattice admitting a nonspecial weakly distinguished basis $B$.
Then the number of orbits of $\Ga_B$ equals $2^\ku+2${\rm,} where 
$\ku=\dim \K$. These orbits are the $2^\ku$ points of $\K$ and the sets 
$q_B^{-1}(0)\setminus \K$ and $q_B^{-1}(1)\setminus \K$.
\endproclaim

\demo{Proof} Evidently, any group generated by transvections acts trivially on
$\K$. Next, in the nonspecial case $q_B^{-1}(1)\setminus \K$ is an orbit
of $\Ga_B$ by Theorem~3.5 of \cite{Ja}. To prove that 
$q_B^{-1}(0)\setminus \K$ is an orbit as well, take an arbitrary pair
$u,v\notin \K$ such that $q_B(u)=q_B(v)=0$. It is easy to see that there exist
$u',v'\notin \K$ such that $q_B(u')=q_B(v')=1$ and $\langle u,u'\rangle=
\langle v,v'\rangle=1$. Define $\Cal D_u=\{w\in \V\: \langle w,u\rangle=
\langle w,u'\rangle=0\}$, $\Cal D_v=\{w\in \V\: \langle w,v\rangle=
\langle w,v'\rangle=0\}$. Evidently, $\K$ is a subspace of both $\Cal D_u$ and
 $\Cal D_v$, and $\dim \Cal D_u=\dim \Cal D_v=d-2$. Let $q_u$ and $q_v$ be 
the restrictions of $q_B$ to $\Cal D_u$ and $\Cal D_v$, respectively. 
If $q_B(\K)=\bF_2$, then the same is 
true for $q_u$ and $q_v$, and so the quadratic spaces $(\Cal D_u,q_u)$ and
$(\Cal D_v,q_v)$ are isometric. Otherwise, if $q_B(\K)=0$, the Arf invariants
for the forms $q_B$, $q_u$, $q_v$ are defined. Moreover, $\Arf(q_u)=\Arf(q_v)=
1+\Arf(q_B)$, and hence  $(\Cal D_u,q_u)$ and $(\Cal D_v,q_v)$ are again 
isometric by the Arf theorem (see \cite{Pf}). In both cases the isometry 
$\Cal D_u\to \Cal D_v$
can be extended to an isometry of the entire $\V$ by letting $u\mapsto v$,
$u'\mapsto v'$. By Theorem~4.1 this isometry belongs to $\Ga_B$. \qed
\enddemo

Let us consider now the conjugate to the action of the monodromy group
of a vanishing lattice. We start from  the following
well-known result.

\proclaim{Lemma~4.4} Let $G$ be a finite group acting on a finite-dimensional
space $\V$ over a finite field $\kk$. Then the number of orbits of the action 
equals the number of orbits of the conjugate action.
\endproclaim

\demo{Proof} The number of orbits of a finite group action on any finite
set $X$ can be calculated using the Frobenius formula (see e.g. \cite{Ke}):
$$
\sharp=\frac{1}{|G|}\sum_{g\in G} |X^g|,
$$ 
where $X^g$ is the set of all $g$-stable points. But the set $\V^g$ of the 
stable points for the linear operator
$E_g$  coincides with $\ker(E_g-E)$, where $E$ is the
identity operator. Therefore, $|\V^g|=|\kk|^{\dim\ker(E_g-E)}$. The equality
$\dim\ker(E_g-E)=\dim\ker(E_g-E)^*= \dim\ker(E^*_g-E)$ proves the statement.
\qed
\enddemo

Consider now the relation between the orbits of a nonspecial monodromy group
$\Ga$ and the orbits of the conjugate action.
Choose some basis $\{b_1,\dots,b_\ku\}$ of $\K$ and a $\ku$-tuple
$\eta=\{\eta_1,\dots,\eta_\ku\}\in\bF_2^\ku$.
Denote by $\A^\eta$ the  affine subspace of $\V^*$ of codimension $\ku$ 
that consists of all elements $x\in \V^*$ such that $(b_j,x)=\eta_j$,
$1\ls j\ls\ku$. Evidently, $\A^\eta$ is invariant under the conjugate
action of $\Ga$ for any $\eta\in \bF_2^\ku$. Moreover, since $\langle
\cdot,\cdot\rangle$ is skew-symmetric, we see that the image 
$L(\V)$ of the map $L\:\V\to \V^*$ coincides with $\A^0$.

Our strategy is as follows. To study the structure of the conjugate action we
use an additional construction. We introduce two subspaces $\V_1\subset \V$ 
and $\V_2\subset \V$, which are transversal to $\K$ (therefore 
$\dim\V_1=\dim\V_2=2m$) and $\V_1+\V_2=\V$. To each subspace $\V_i$
corresponds a subgroup $\Ga_i\subset \Ga$ (generated by transvections w.r.t.\
the elements in $\V_i$). One can easily see that the set
$\Ga_1\cup\Ga_2$ generates $\Ga$. We first study the orbits of $\Ga_i$
and of the corresponding conjugate action separately (which is fairly simple),
and then describe their interaction. We start with the case of one subspace.

Consider a subspace $\V_1\subset \V$ transversal to $\K$. The restriction of
the form $\langle\cdot,\cdot\rangle$ is nondegenerate on $\V_1$. Let
$\Gamma_1\subset \Gamma$ be the subgroup generated by the tranvections w.r.t.
elements of $\V_1$.

\proclaim{Lemma 4.5} Let $\Ga$ be nonspecial. Then
for any $\eta\in\bF_2^\ku$ the affine subspace 
$\A^\eta$ consists of three orbits of the conjugate $\Ga_1$-action.
\endproclaim

\demo{Proof} Observe first that for any $\eta\in\bF_2^\ku$
the intersection of $\A^\eta$ and $\V_1^\perp \subset\V^*$ contains 
exactly one element, which we denote by $A_1^\eta$. Indeed, if
$\{v_1,\dots,v_{2m}\}$ is an arbitrary basis of $\V_1$, then $A_1^\eta$ is
the unique solution of the equations $(b_j,A_1^\eta)=\eta_j$, $1\ls
j\ls\ku$, $(v_j,A_1^\eta)=0$, $1\ls j\ls 2m$.

Next, let $L_1$ be the restriction of $L$ to $\V_1$, and $L_1^\eta$ be the 
composition of $L_1$ with the translation by $A_1^\eta$. We say that
$A_1^\eta$ is the {\it shift\/} corresponding to $\A^\eta$; observe that
$A_1^0=0$, and thus $L_1^0=L_1$. Evidently, 
$L_1^\eta$ provides an affine isomorphism between $\V_1$ and $\A^\eta$.
Recall that the conjugate action of $\Gamma$ (and hence, of $\Gamma_1$)
preserves $\A^\eta$. Moreover, since $\Gamma_1$ preserves $\V_1$, the
element $A_1^\eta$ is a fixed point of the conjugate to the
$\Ga_1$-action, as the unique annihilator of $\V_1$ in $\A^\eta$.
Therefore, the diagram
$$
\CD
\V_1 @>{L_1^\eta}>> \A^\eta \\
@V g_1 VV @A g_1^* AA \\
\V_1 @>{L_1^\eta}>> \A^\eta
\endCD
$$
is commutative for any $g_1\in \Gamma_1$. Since $\Ga$ is nonspecial, the
same is true for $\Ga_1$ (by Lemma~4.2). It remains to apply Lemma~4.3
with $\ku=0$.
\qed
\enddemo

Let us now choose one more subspace $\V_2\subset \V$ transversal to $\K$
such that $\V_1+\V_2=\V$ and $\Ga_1\cup\Ga_2$ generates $\Ga$. To study 
the orbits of the conjugate $\Ga$-action we have to find which orbits
of the conjugate $\Ga_1$-action are glued together by the conjugate
$\Ga_2$-action. Let us define an affine isomorphism $I^\eta\: \V_1\to
\V_2$ by $I^\eta=(L_2^\eta)^{-1}\circ L_1^\eta$, where $L_2$ is the
restriction of $L$ to $\V_2$.

\proclaim{Lemma~4.6} Let $\Ga$ be nonspecial.

{\rm (\romannumeral1)} For $\eta=0$ the following 
alternative holds\/{\rm:} the linear space $\A^0$ consists of two 
orbits of the conjugate $\Gamma$-action if there exist $u,v\in \V_1$
such that $q_B(u)=0${\rm,} $q_B(v)=q_B(I^0(u))=q_B(I^0(v))=1${\rm,}
and of three orbits of the conjugate $\Gamma$-action otherwise.

{\rm (\romannumeral2)} 
For any $\eta\ne0$ the following alternative holds\/{\rm:}
the affine subspace $\A^\eta$ consists of one orbit of the conjugate
$\Gamma$-action if there exist $u,v\in \V_1$
such that $q_B(u)=0${\rm,} $q_B(v)=q_B(I^\eta(u))=q_B(I^\eta(v))=1${\rm,}
and of two orbits of the conjugate $\Gamma$-action otherwise.
\endproclaim

\demo{Proof} The first part of this lemma follows easily from Lemma~4.3
with $\ku=0$ together with Lemma~4.5.
To prove the second part we notice first that $I^\eta(0)\ne 0$, 
and then we are done by the same reasons.
\qed
\enddemo

\heading \S 5. Orbits of the second $\G_n$-action \endheading

In this section we prove Theorem~2.5 describing the orbits of
the second $\G_n$-action on the space of upper-triangular $(n-1)\times (n-1)$ 
matrices over $\bF_2$. Observe that the basic convention in this section is 
different from that of \S 4, namely, the conjugate to the second $\G_n$-action
is the action  preserving a natural skew-symmetric form, and we study 
the second $\G_n$-action itself as the conjugate to its conjugate using
Lemmas~4.4-4.6.

\subheading {5.1. The conjugate to the second $\G_n$-action} 
In this section we assume that $\V=(T^{n-1}(\bF_2))^*$. 
Let us introduce a bilinear form $\langle\cdot,\cdot\rangle$ on $\V$ by
$\langle M, N\rangle=\sum_{(i_1,j_1),(i_2,j_2)} m_{i_1,j_1}n_{i_2,j_2}$,
where $(i_1,j_1),(i_2,j_2)$ runs over all pairs of neighbors in $\H_{n-1}$,
see \S 3. Evidently, $\langle\cdot,\cdot\rangle$ is skew-symmetric.

Let $B_{n-1}$ be the standard basis of $\V$ consisting of matrix entries.
Elements of $B_{n-1}$ correspond bijectively to the vertices of $\H_{n-1}$.
Denote by $Q$ the quadratic function $q_{B_{n-1}}$ associated with 
$\langle\cdot,\cdot\rangle$. To find the value $Q(M)$ for an arbitrary 
matrix $M\in \V$ we represent $M$ as $M=\sum_{b_i\in B_{n-1}}\alpha_ib_i$
and define $\Supp M=\{b_i\in B_{n-1}\: \alpha_i=1\}$. Let $\gr(M)$ denote
the subgraph of $\H_{n-1}$ induced by the vertices corresponding to
$\Supp M$. The following statement follows easily from the definitions.

\proclaim{Lemma~5.1} The value $Q(M)$ equals $\mod 2$ the number of vertices
of $\gr(M)$ plus the number of its edges.
\endproclaim

As a corollary, we obtain the values of $Q$ on the dual invariants 
$\widetilde P_1,\dots,\widetilde P_k$.

\proclaim{Corollary~5.2} {\rm (\romannumeral1)} Let $n=2k+1${\rm,} 
then $Q(\widetilde P_i)=0$ for $1\ls i\ls k$.

{\rm (\romannumeral2)} Let $n=2k${\rm,} then $Q(\widetilde P_{i})=1$ for 
$1\ls i\ls k-1${\rm,} and $Q(\widetilde P_k)=k \mod 2$.
\endproclaim

\demo{Proof} (i) It follows from the description
of the patterns $\pi_i$ that $\gr(P_i)$, $1\ls i<k$, consists
of three copies of $\H_i$, several disjoint cycles, and one more copy of
$\H_l$ with $l=|n-3i-2|+1$, provided $j_i=\min\{i+1,n-2i-1\}$ is even.
Therefore, by Lemma~5.1, we have $Q(P_i)=3i(2i-1) \mod 2$ for $j_i=i+1$, 
$i$ even, and $Q(P_i)=3i(2i-1)+l(2l-1) \mod 2$ otherwise. Since $n$ is odd, 
we have $l=i \mod 2$, and hence $Q(P_i)=0$, $1\ls i<k$. Finally, for $P_k$ 
we have  $Q(P_k)=(n-1)(2n-3)=2k(4k-1)\mod 2=0$. It remains to notice that
$Q(\widetilde P_i)=Q(P_i)+Q(P_{i-1})$.

(ii) In the same way as in (i) we see that if $1\ls i<k$,then
$Q(P_i)=3i(2i-1)+l(2l-1) \mod 2$ for $j_i=i+1$, $i$ odd, and 
$Q(P_i)=3i(2i-1) \mod 2$ otherwise. Since $n$ is even, we have 
$l=i+1 \mod 2$, and hence $Q(P_i)=0$ for even $i<k$ and $Q(P_i)=1$ 
for odd $i<k$. Finally, for $P_k$ we have $Q(P_k)=(n-1)(2n-3)=
(2k-1)(4k-3)\mod 2=1$. As in (i), the result follows from 
$Q(\widetilde P_i)=Q(P_i)+Q(P_{i-1})$.
\qed
\enddemo

Let $\K$ be the kernel of $\langle\cdot,\cdot\rangle$. Evidently, $\K$ is
the space of dual invariants of the second $\G_n$-action, hence Theorem~2.3
gives an explicit description of $\K$. Applying Corollary~5.2 we get the 
following statement.

\proclaim{Corollary~5.3} Let $n=2k${\rm,} then $Q(\K)=\bF_2$ and
$|Q^{-1}(0)|=|Q^{-1}(1)|=2^{2k^2-k-1}${\rm,}
$|Q^{-1}(0)\cap\K|=|Q^{-1}(1)\cap\K|=2^{k-1}$.   
\endproclaim
   
\demo{Proof}  It follows immediately from Corollary~5.2(ii) that $Q(\K)=\bF_2$.
The statement concerning the sizes of $|Q^{-1}(0)|$ and $|Q^{-1}(1)|$ 
follows easily from the general description (see \S 4.1)
with $m=k(k-1)$ and $\ku=k$. The last statement follows from the fact that
the normal form of $Q|_\K$ in case $Q(\K)=\bF_2$ is $x_1^2$ (see \cite{Pf}).
\qed
\enddemo

In case $n=2k+1$ the situation is far more complicated.

\proclaim{Lemma~5.4} Let $n=2k+1${\rm,} then $Q(\K)=0$. If $k=4t+1${\rm,}
then $\Arf(Q)=1$ and $|Q^{-1}(0)|=2^{2k^2+k-1}-2^{k^2+k-1}${\rm,}
$|Q^{-1}(1)|=2^{2k^2+k-1}+2^{k^2+k-1}$. Otherwise $\Arf(Q)=0$ and 
$|Q^{-1}(0)|=2^{2k^2+k-1}+2^{k^2+k-1}${\rm,}
$|Q^{-1}(1)|=2^{2k^2+k-1}-2^{k^2+k-1}$.
\endproclaim
   
\demo{Proof} It follows immediately from Corollary~5.2(i) that $Q(\K)=0$,
and hence the Arf invariant exists. Observe that by general theory of 
quadratic spaces (see \S4.1) this means that $|Q^{-1}(0)|\ne|Q^{-1}(1)|$.
Moreover, $|Q^{-1}(0)|>|Q^{-1}(1)|$ implies $\Arf(Q)=0$,
while $|Q^{-1}(0)|<|Q^{-1}(1)|$ implies $\Arf(Q)=1$. Therefore, to find 
$\Arf(Q)$ it suffices to count the number of elements in $\V$ on which
$Q$ vanishes.

Let $\omega\:\V\to\bF_2^{n-1}$ denote the projection on the first row (see 
Lemma~3.3), and let $\V_a=\omega^{-1}(a)$, $a\in\bF_2^{n-1}$. We say that
$\V_a$ is {\it inessential\/} if $|Q^{-1}(0)\cap\V_a|=|Q^{-1}(1)\cap\V_a|$,
and {\it essential\/} otherwise.

Let us choose 
$M_a\in\V_a$ such that $Q(M_a)=0$; if $\gr(a)$ contains an even number of 
connected components, one can take $\Supp M_a=\Supp a$, otherwise 
$\Supp M_a$ contains one more vertex, which corresponds to the matrix entry 
$(n,n)$. We thus have $M_a=M_a^0+M_a^1$, where $\Supp M_a^0 =\Supp a$ and
$|\Supp M_a^1|\ls1$. 

Let us define a function $Q_a$ on $\V_0$ by $Q_a(M)=Q(M+M_a)$ for any 
$M\in\V_0$; observe that $Q_a$ is just a shift of the restriction $Q|_{\V_a}$.
Evidently, $Q_a(M)=Q(M)+\langle M,M_a\rangle$. Therefore, $Q_a(M+N)-Q_a(M)-
Q_a(N)=\langle M,N\rangle$, which means that $Q_a$ and $Q$ define the same 
bilinear form on $\V_0$ (observe that $\V_0$ is identified naturally with 
$(T^{n-2}(\bF_2))^*$). 
Let us evaluate $Q_a$ on the kernel $\K_0$; $Q_a(\K_0)=\bF_2$ would
mean that $Q$ vanishes exactly on a half of the elements of $\V_a$, in 
other words, that $\V_a$ is inessential, while $Q_a(\K_0)=0$ would mean that $\V_a$
is essential. 

Since both $Q$ and $Q_a$ define the same bilinear form on $\V_0$, we see
that $\K_0$ is the space of dual invariants of the second $\G_{n-1}$-action;
therefore, by Theorem~2.3, it has a basis $\{\widetilde P_1,\dots,
\widetilde P_k\}$. By Corollary~5.2(ii),
$Q_a(\widetilde P_i)=1+\langle \widetilde P_i,M_a\rangle=1+
\langle \widetilde P_i,M_a^0\rangle$ for $1\ls i\ls k-1$ and 
$Q_a(\widetilde P_k)=k+\langle \widetilde P_k,M_a^0\rangle$.
Therefore, $\V_a$ is essential if and only if the entries of $a$ satisfy
equations $a_{k-i}+a_{k+i+1}=k+i\mod 2$ for $0\ls i\ls k-1$.
It is easy to see that any solution of the above equations is represented
as $a=\bar h+s$, where $\bar h$ is as defined in Theorem~2.2(ii) and 
$s\in\Sym^{n-1}$. It follows from Theorem~2.3 that for any $s\in\Sym^{n-1}$
there exists $S\in\K$ such that $\omega(S)=s$. Observe that in this case
$M\mapsto M+S$ takes $\V_a$ to $\V_{a+s}$. Moreover, $Q(M+M_a+S)=Q(M+M_a)$,
since $\langle M+M_a,S\rangle=0$ follows from $S\in\K$ and $Q(S)=0$ by
Corollary~5.2(i). Therefore, all the essential subspaces $\V_a$ influence 
$\Arf(Q)$ in the same way, hence, $\Arf(Q;\V)=\Arf(Q_{\bar h};\V_0)$. 

To study $\Arf(Q_{\bar h})$ on $\V_0$ we reiterate the same process once more,
that is, we decompose $\V_0$ into affine subspaces $\V_{0b}$, where $b$ is
defined by the projection $\omega_0$ of $\V_0$ on the first row. We choose 
$M_{0b}\in\V_{0b}$ similarly to $M_a$ and define a function
$Q_{0b}$ on $\V_{00}$ by $Q_{0b}(M)=Q_{\bar h}(M+M_{0b})$. As before, 
$\V_{00}$ is identified naturally with $(T^{n-3}(\bF_2))^*$, and 
$Q_{0b}$ and $Q$ define the same bilinear form on $\V_{00}$. 
The corresponding kernel $\K_{00}$ is the space of dual 
invariants of the second $\G_{n-2}$-action; its basis is 
$\{\widetilde P_1,\dots,\widetilde P_{k-1}\}$. We thus get that $\V_{0b}$
is essential if and only if the entries of $b$ satisfy equations $b_{k-i}+
b_{k+i}=0$ for $1\ls i\ls k-1$. Any solution of these equations belongs to
$\Sym^{n-2}$. As before, $M\mapsto M+S$ with $S\in\K_0$ and $\omega_0(S)=s$
takes $\V_{0b}$ to $\V_{0,b+s}$ and $Q_{\bar h}(M+M_{0b}+S)=
Q_{\bar h}(M+M_{0b})$; identity $Q_{\bar h}(S)=0$ follows readily from
$Q(\widetilde P_i)=\langle\widetilde P_i,M_{\bar h}^0\rangle$, $1\ls i\ls k-1$.
We thus get $\Arf(Q_{\bar h};\V_0)=\Arf(Q_{\bar h};\V_{00})$.

Recall that on $\V_{00}$ one has $Q_{\bar h}(M)=Q(M)+
\langle M,M_{\bar h}^1\rangle$. It is easy to see that for $k=4t$ and $k=4t+3$
one has $M_{\bar h}^1=0$, and hence $Q_{\bar h}\equiv Q$ on $\V_{00}$. This
means that the Arf invariant is constant on each  triple of the form $k=4t+2,
4t+3,4t+4$. 

Let now $k=4t+1$ or $k=4t+2$, and hence $M_{\bar h}^1\ne0$.
We decompose $\V_{00}$ into four affine subspaces $\V_{00}^{00}$, 
$\V_{00}^{01}$, $\V_{00}^{10}$, and $\V_{00}^{11}$. The subspace $\V_{00}^{ij}$
consists of the matrices having $i$ at position $(n-1,n-1)$ and $j$ at
position $(n-1,n)$. It is easy to see that $Q_{\bar h}\equiv Q$ on 
$\V_{00}^{ii}$; moreover, the involution $M\mapsto M+M_{\bar h}^1$ reverses
the value of $Q$ (and thus of $Q_{\bar h}$) on these subspaces. Therefore, 
the subspaces $\V_{00}^{ii}$ are inessential. On the other hand,  
$Q_{\bar h}\equiv Q+1$  on $\V_{00}^{ij}$, $i\ne j$, and  the involution 
$M\mapsto M+M_{\bar h}^1$ in this case preserves the value of $Q$ 
(and thus of $Q_{\bar h}$). Therefore, $\Arf(Q_{\bar h};\V_{00})=1+
\Arf(Q;\V_{00})$. This means that the Arf invariant reverses twice on each
triple of the form $k=4t,4t+1,4t+2$.

To complete the proof it is enough to check the value of the Arf invariant
for $k=1$.

The exact values for the lengths of orbits follow easily from the general 
description (see \S 4.1) with $m=k^2$ and $\ku=k$.
\qed
\enddemo

Let us now find the relation between the conjugate to the second 
$\G_n$-action and the theory of skew-symmetric vanishing lattices explained
in \S 4.

\proclaim{Lemma~5.5} {\rm (\romannumeral1)} The triple $(\V,
\langle\cdot,\cdot\rangle, B_{n-1})$ is a vanishing lattice.

{\rm (\romannumeral2)}
The conjugate to the second $\G_n$-action coincides with the action
of $\Ga_{B_{n-1}}$.

{\rm (\romannumeral3)} The basis $B_{n-1}$ is weakly distinguished and its
graph $\gr(B_{n-1})$ coincides with $\H_{n-1}$.
\endproclaim

\demo{Proof} (i) We have to check that $B_{n-1}$ satisfies conditions (i)--(iii)
of the definition of the vanishing lattices. Condition (ii) is evident.
To check to it suffices to take for $\delta_1$ and $\delta_2$ any pair of
adjacent vertices of $\H_{n-1}$. Finally, to check (i) we take an arbitrary 
pair $b$, $b'$ of adjacent vertices of $\H_{n-1}$ and find that 
$T_bT_{b'}(b)=b'$.

(ii) Follows immediately from Lemma~3.2.

(iii) Obvious.
\qed\enddemo

\proclaim{Lemma~5.6} The conjugate to the second $\G_n$-action is the action
of a nonspecial monodromy group for $n\gs 5$.
\endproclaim

\demo{Proof} By Lemmas~4.2 and~5.5 we have to check that the Dynkin diagram
of $E_6$ is an induced subgraph of $\H_m$ for $m\gs4$. Since $\H_m$ is an 
induced subgraph of $\H_l$ for $m<l$, it suffices to find an induced subgraph 
corresponding to $E_6$ in $\H_4$, see Fig.~4.
\qed\enddemo

\vskip 10pt
\centerline{\hbox{\hskip 0.2cm \epsfysize=3cm\epsfbox{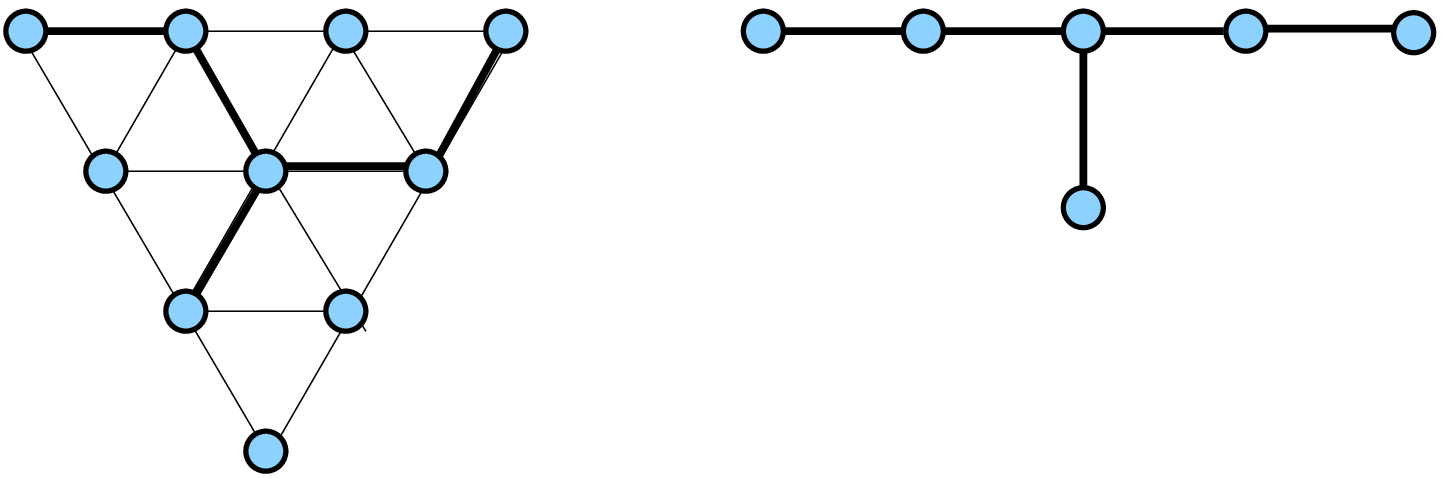}}}
\midspace{1mm}
\caption{Fig.~4. The graph $\H_4$ and the Dynkin diagram for $E_6$}
\vskip 5pt

\subheading{5.2. The orbits of the second $\G_n$-action}
First of all, let us find the number of orbits of the second $\G_n$-action.

\proclaim{Lemma~5.7} The number of orbits of the second $\G_n$-action equals
$2^k+2$ for $n\gs5$.
\endproclaim

\demo{Proof} Indeed, by Lemma~4.4 the number of orbits of the second 
$\G_n$-action equals that of its conjugate. By Lemma~5.6 the conjugate 
to the second $\G_n$-action is the action of a nonspecial monodromy
group; hence, by Lemma~4.3 the number of its orbits equals
$2^\ku+2$, where $\ku$ is the dimension of $\K$. Since $\K$ is just
the space of dual invariants of the second $\G_n$-action, 
Theorem~2.3(ii) implies $\ku=k$.
\qed
\enddemo

Let us now define two subsets of the vertex set of $\H_{n-1}$ as follows:
$\Sigma_1$ consists of the vertices corresponding to the matrix entries 
$(i,n-1)$, $1\ls i\ls k$, and $\Sigma_2$ of the vertices corresponding
to $(n-i,n-i)$, $1\ls i\ls k$, see Fig.~5.  

\vskip 10pt
\centerline{\hbox{\hskip 0.2cm \epsfysize=3cm\epsfbox{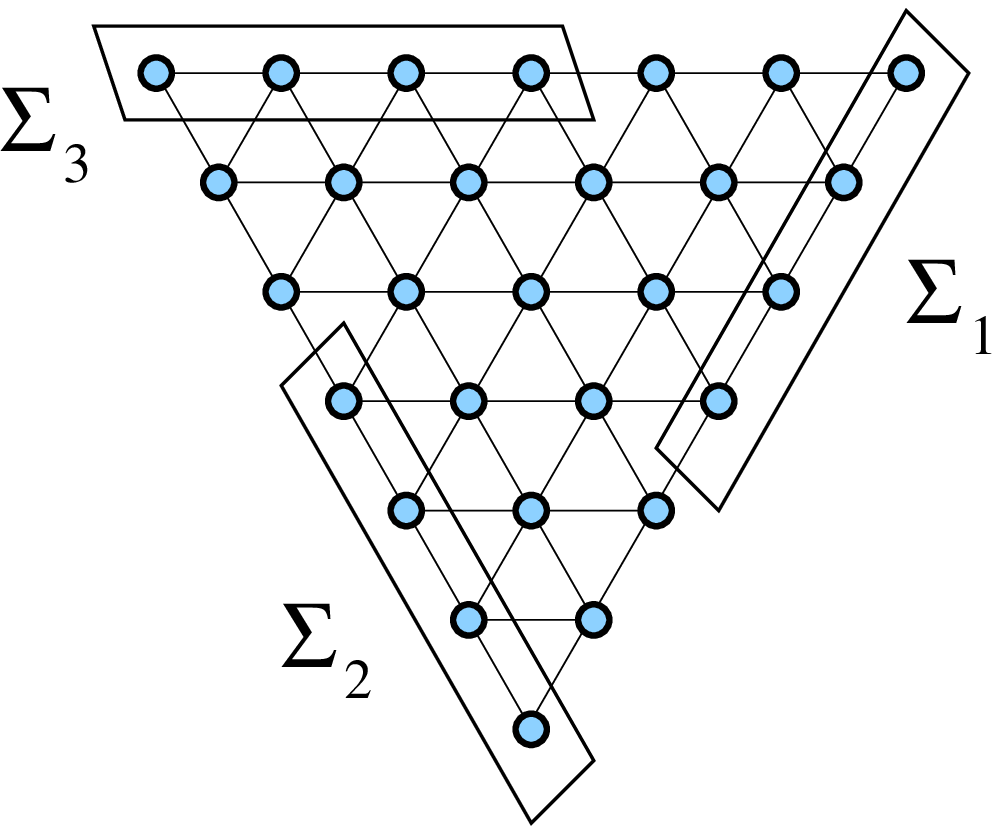}}}
\midspace{1mm} \caption{Fig.~5. The sets $\Sigma_1$, $\Sigma_2$, $\Sigma_3$ 
for the case $n=8$}
\vskip 5pt

Following the strategy described in \S 4, we choose the following two subspaces
in $\V=(T^{n-1}(\bF_2))^*$:
$$\gather
\V_1=\{M\in \V\: \Supp M\cap \Sigma_1=\varnothing\},\\
\V_2=\{M\in \V\: \Supp M\cap \Sigma_2=\varnothing\}.\endgather
$$
Let $\K$ denote the kernel of $\langle\cdot,\cdot\rangle$.
It follows easily from the explicit description of $\K$ (see Theorem~2.3(ii))
that $\codim \V_1=\codim \V_2=\dim \K$ and $\V_1\cap \K=\V_2\cap \K=0$, hence
$\V_1$ and $\V_2$ as above satisfy the assumptions of \S 4.

Fix the basis $\widetilde P_1,...,\widetilde P_k$ of $\K$.
It is easy to see that the stratum of the second $\G_n$-action at height
$\eta$ is just the affine subspace $\A^\eta$ as defined in \S 4.2.
Therefore, the statements of Theorem~2.5(i) concerning the number and the
lengths of orbits are yielded by the following proposition.

\proclaim{Lemma~5.8} Let $n=2k+1\gs5$. 

{\rm (\romannumeral1)} The linear subspace $\A^0\subset \V^*$ consists of
three orbits of the second $\G_n$-action. The lengths of the orbits are 
$2^{2k^2-1}-\ee_k2^{k^2-1}${\rm,} $2^{2k^2-1}+\ee_k2^{k^2-1}-1${\rm,}
and $1${\rm,} where $\ee_k=-1$ for $k=4t+1$ and $\ee_k=1$ otherwise.

{\rm (\romannumeral2)} For any nonzero $\eta\in\bF_2^{k}$
the affine subspace $\A^\eta\subset \V^*$ is an orbit of 
the second $\G_n$-action.
\endproclaim

\demo{Proof} (i) Let us first find the number of orbits contained in $\A^0$.
By Lemma~4.6(i) it suffices to prove that
$Q(M)=Q(I^0(M))$ for any $M\in \V_1$. Define $M^\K=M+I^0(M)$, and 
let $L\: \V\to \V^*$ be the linear mapping associated with $\langle\cdot,\cdot
\rangle$. Since by definition $L(M)=L(I^0(M))$ for any $M\in \V_1$, we
have $L(M^\K)=0$, and thus $M^\K\in \K$. Therefore $Q(I^0(M))=Q(M)+Q(M^\K)$,
and since by Lemma~5.4 $Q(\K)=0$, we are done.

It follows now from Lemmas~4.5 and~4.3 that the three orbits in question are
isomorphic to $\{0\}$, $(Q^{-1}(0)\cap \V_1)\setminus\{0\}$ and 
$Q^{-1}(1)\cap \V_1$. Evidently, $\Arf(Q)=\Arf(Q|_{\V_1})$ and 
$(\V_1,Q|_{\V_1})$ is a nondegenerate quadratic space of dimension $2k^2$.
Therefore the statement concerning the lengths of the orbits follows from
Lemma~5.4 and the general description of quadratic spaces (see \S 4.1) 
with $m=k^2$ and $\ku=0$.

(ii) By Lemma~5.7 the total number of orbits is $2^k+2$, and by 
part (i) of the present Lemma exactly three orbits are contained in $\A^0$. 
Therefore, each of the $2^k-1$ remaining orbits coincides with the 
corresponding $\A^\eta$.
\qed
\enddemo

The corresponding statements of Theorem~2.5(ii) are yielded by the following 
proposition.

\proclaim{Lemma~5.9} Let $n=2k\gs6$  and 
$\bar\eta=(Q(\widetilde P_1),\dots,Q(\widetilde P_k))$.
 
{\rm (\romannumeral1)} The affine subspace $\A^{\bar\eta}\subset \V^*$ 
consists of two orbits of the second $\G_n$-action. The lengths of the orbits 
are $2^{2k(k-1)-1}-2^{k(k-1)-1}$ and $2^{2k(k-1)-1}+2^{k(k-1)-1}$.

{\rm (\romannumeral2)} The linear subspace $\A^0\subset \V^*$ consists of
two orbits of the second $\G_n$-action. The lengths of the orbits are 
$2^{2k(k-1)}-1$ and $1$.

{\rm (\romannumeral3)} For any $\eta\in\bF_2^{k}${\rm,} 
$\eta\ne0,\bar\eta${\rm,} the affine subspace $\A^\eta\subset \V^*$ is an
orbit of the second $\G_n$-action.
\endproclaim

\demo{Proof} (i) Let $A_1^\eta\in \V^*$ be the shift corresponding to
$\V_1$, as defined in the proof of Lemma~4.5, and $A_2^\eta$ be the
similar shift corresponding to $\V_2$. 
 Observe that $A_1^\eta+A_2^\eta\in \A^0$, and hence there exists 
$H_2^\eta\in \V_2$ such that $L(H_2^\eta)=A_1^\eta+A_2^\eta$. 

For any $M\in \V_1$ put $M^\K=M+I^\eta(M)+H_2^\eta$. It follows from the 
definition of $I^\eta$ that $L(I^\eta(M))=L(M)+A_1^\eta+A_2^\eta$, hence
$L(M^\K)=0$, that is, $M^\K\in \K$. 

Let us now find $Q(I^{\bar\eta}(M))$. Since $I^\eta(M)=M+H_2^\eta+M^\K$
and $\langle M,H_2^\eta\rangle=(M,A_1^\eta+A_2^\eta)$,
one has $Q(I^\eta(M))=Q(M)+Q(H_2^\eta)+Q(M^\K)+(M,A_1^\eta+A_2^\eta)$,
where $(\cdot,\cdot)$ is the standard coupling between $\V$ and $\V^*$
defined in \S 2.

Since $M\in\V_1$, we get $(M,A_1^\eta+A_2^\eta)=(M,A^\eta_2)=
|\Supp M\cap\Supp A_2^\eta|$. It follows easily from Theorem~2.3(ii) that
$\Supp A_2^{\bar\eta}=\{(n-i,n-i)\in \Sigma_2\: Q(\widetilde P_i)=1\}$.
We thus have $(M,A_1^{\bar\eta}+A_2^{\bar\eta})=|\{(n-i,n-i)\in\Sigma_2
\cap\Supp M\: Q(\widetilde P_i)=1\}|$.

On the other hand, 
$M^\K=\sum\{\widetilde P_i\: (n-i,n-i)\in\Sigma_2\cap\Supp M^\K\}$, and
hence
$Q(M^\K)=|\{(n-i,n-i)\in\Sigma_2\cap\Supp M^\K\: Q(\widetilde P_i)=1\}|$. 
Finally, $M+M^\K= I^\eta(M)+H_2^\eta\in\V_2$, and hence 
$\Sigma_2\cap\Supp M=\Sigma_2\cap\Supp M^\K$.

Thus, $Q(I^{\bar\eta}(M))-Q(M)=Q(H_2^{\bar\eta})$ does not depend on $M$,
and hence $u$,
$v$ as in Lemma~4.6(ii) do not exist. Therefore, $\A^{\bar\eta}$ consists
of the two orbits of the second $\G_n$-action. For the same reason the trivial
one-element orbit obtained from $\{0\}$ is glued to the image of 
$Q^{-1}(0)\setminus\{0\}$, and hence the lengths of the orbits
are the sizes of $Q^{-1}(0)$ and $Q^{-1}(1)$ for the nondegenerate 
case. Thus $\ku=0$ and $m=k(k-1)$, which yield the required result.

(ii), (iii) Follows immediately from Lemma~5.7 and part (i) in the same way as
 in Lemma~5.8(ii). To find the lengths of the orbits in $\A^0$ observe that
$L_1(0)=L_2(0)=0\in\A^0$, hence the two parts that are glued together are
$Q^{-1}(0)\setminus\{0\}$ and $Q^{-1}(1)$, and the result follows.
\qed
\enddemo 

To obtain Theorem~2.5(ii) it suffices to notice that by Corollary~5.2(ii)
$\bar\eta$ as defined in Lemma~5.9 coincides with $\bar\eta$ as defined
in Theorem~2.5(ii). 

\heading \S 6. Orbits of the first $\G_n$-action \endheading

In this section we prove Theorem~2.2 concerning the structure of 
orbits of the first $G_n$-action.

\subheading{6.1. The structure of the strata of the first $\G_n$-action}
Let us introduce some notation. We write $X\sim\widetilde X$ if the
matrices $X$ and $\widetilde X$ belong to the same orbit. We do not
specify in this notation which $\G_n$-action is meant; this should be 
always clear from the context, since each space under consideration
carries exactly one $\G_n$-action. 

Let $\Omega=\omega_1\omega_2\dots\omega_t$ denote an arbitrary word in
alphabet $\{g_{ij}\}$; we write $\bar\Omega$ for the word 
$\omega_t\dots\omega_2\omega_1$. For any subset $\Sigma$ of the alphabet 
we define $[\Omega:\Sigma]$ as the number of occurencies of $g_\sigma$,
$\sigma\in\Sigma$, in $\Omega$. We write $\Omega X$ to denote the matrix
obtained from $X$ by applying to it the elements $\omega_1,\omega_2,\dots,
\omega_t$ in this order. We say that $\Omega$ is $X$-nonredundant if
$\{\omega_1\dots\omega_s\}X\ne\{\omega_1\dots\omega_{s-1}\}X$ for
$s=1,\dots,t$ (here $\omega_0X=X$). Evidently, $X$-nonredundancy is
preserved under isomorphisms and under $\Psi_n$.

It is easy to see that each stratum $S^h$ of the first $\G_n$-action is a
covering of degree $2^k$ over the stratum $\Psi_n(S^h)$ of the second
$\G_n$-action. We are going to find a family of linear functionals $C^h\:
T^n(\bF_2)\to\bF_2$ with the following property: let $M,\widetilde M\in
S^h$ be an arbitrary pair of matrices such that
$\Psi_n(M)=\Psi_n(\widetilde M)$, then $M\sim\widetilde M$ 
if and only if $C^h(M)=C^h(\widetilde M)$. The existence of such a family
would imply immediately that each nontrivial orbit of the first
$\G_n$-action is a covering over the corresponding orbit of the second
$\G_n$-action of degree $2^{k-1}$, except for the cases when $C^h$ is
trivial; in the latter case the degree equals $2^k$.

For any $M\in T^n(\bF_2)$ we denote by $S(M)$ the set of matrices
$\widetilde M$ such that $\Psi_n(M)=\Psi_n(\widetilde M)$ and
$M,\widetilde M$ belong to the same stratum of the first $\G_n$-action.
Further, let $\W$ denote the $k$-dimensional subspace of $T^n(\bF_2)$
generated by the entries $(1,i)$, $1\ls i\ls k$, and $\tau\:T^n(\bF_2)\to
\W$ denote the natural projection.

\proclaim{Lemma~6.1} The projection $\tau$ provides an affine 
isomorphism between $S(M)$ and $\W$ for any $M\in T^n(\bF_2)$.
\endproclaim

\demo{Proof} Indeed,
$S(M)=M+(\ker\Psi_n\cap\DD_n^\perp)=M+(\Cal I_n\cap\DD_n^\perp)$. It follows
from the explicit description of $\DD_n$ (see Theorem~2.1(ii)) that 
$\dim(\Cal I_n\cap\DD_n^\perp)=k$ and $\ker\tau|_{\Cal I_n\cap\DD_n^\perp}=0$.
\qed
\enddemo

Let us introduce, in addition to $\Sigma_1$ and $\Sigma_2$, one more subset 
of the vertex set of $\H_{n-1}$: $\Sigma_3$, consisting of vertices 
corresponding to matrix entries $(1,i)$, $1\ls i\ls k$, see Fig.~5.
As in \S 5.2, we define the subspace $\V_3\subset\V=(T^{n-1}(\bF_2))^*$,
$$
\V_3=\{M\in\V\: \Supp M\cap\Sigma_3=\varnothing\}.
$$
Recall that $L$ provides an isomorphism between $\V_3$ and 
$\A^0=\D_{n-1}^\perp$, and $A_1^\eta+A_2^\eta\in\A^0$; therefore, there
exists $H^\eta_3\in\V_3$ such that $L(H^\eta_3)=A^\eta_1+A^\eta_2$.

\proclaim{Lemma~6.2} $Q(H^\eta_3)+(A^\eta_2,H^\eta_3)=0$ for any
$\eta\in\bF_2^k$.
\endproclaim

\demo{Proof} Let us give an explicit description of $H^\eta_3$. We define
the matrices $H_i\in\V$, $1\ls i\ls k$, in the following way. Let $1\ls i
<k$, then $\Supp H_i$ is obtained from the initial shape by deleting 
the upper-left justified $(n-i-1)\times(n-i-1)$ triangle and upper-right
justified $(i-1)\times(i-1)$ triangle. For $n=2k+1$ we define $H_k$ in the
same way, and for $n=2k$ $\Supp H_k$ is the lower-right justified $(k-1)\times
(k-1)$ triangle, see Fig.~6.

\vskip 20pt
\centerline{\hbox{\hskip 0.2cm \epsfysize=3cm\epsfbox{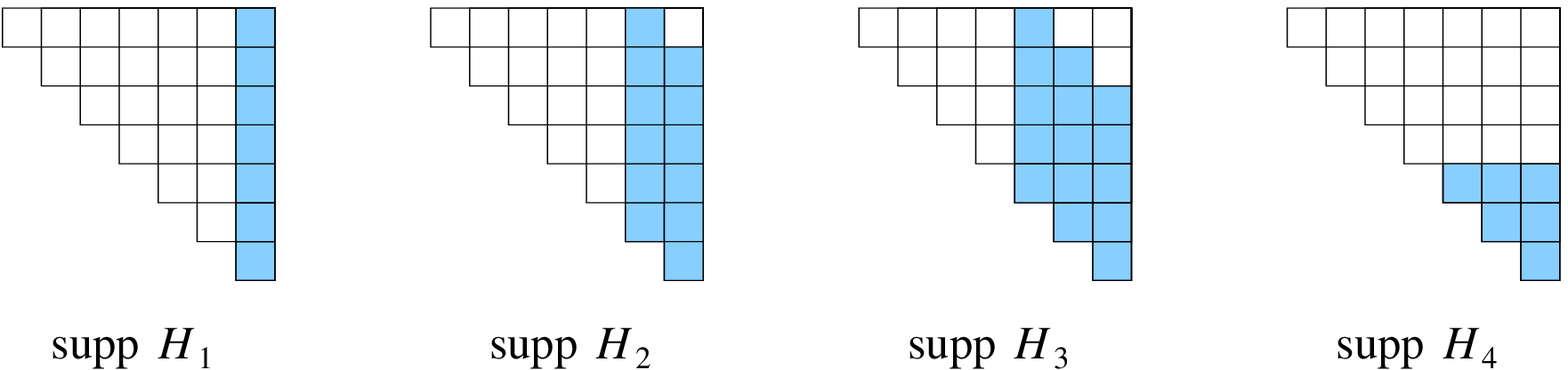}}}
\midspace{1mm} \caption{Fig.6. Supports of the matrices $H_i$ in case $n=8$}
\vskip 5pt

Evidently, $H_i\in\V_3$, $1\ls i\ls k$. Besides, $Q(H_i)=1$ for $1\ls i <k$
and $Q(H_k)=n \mod 2$. Finally, $\langle H_i,H_j\rangle=1$ for $i\ne j$.

Observe now that $H^\eta_3=\sum\{H_i\: (n-i,n-i)\in\Supp A^\eta_2\}$.
Therefore, for $n=2k+1$ we have $Q(H^\eta_3)=c+c(c-1)/2 =c(c+1)/2 \mod 2$, 
where $c=|\Supp A^\eta_2|$. On the other hand, 
$(A^\eta_2,H^\eta_3)=1+2+\cdots+c=c(c+1)/2 \mod 2$. 
Hence $Q(H^\eta_3)+(A^\eta_2,H^\eta_3)=c(c+1)\mod 2=0$.

Let now $n=2k$. If $(n-k,n-k)\notin \Supp A^\eta_2$, then the above proof
applies. Otherwise $Q(H^\eta_3)=c-1+c(c-1)/2=(c-1)(c+2)/2 \mod 2$ and
$(A^\eta_2,H^\eta_3)=1+2+\cdots+c-1+c-1=(c-1)(c+2)/2 \mod 2$, and hence
$Q(H^\eta_3)+(A^\eta_2,H^\eta_3)=(c-1)(c+2)\mod 2=0$.
\qed
\enddemo

Let us define a family of functionals $f^\eta\:\K\to\bF_2$, $\eta\in\bF_2^k$,
as follows: $f^\eta(X)=Q(X)+(A_2^\eta,X)$. It is easy to see that $f^\eta$ is
a linear functional for any $\eta$, since $Q$ is linear on $\K$. The following
statement stems immediately from Corollary~5.2.

\proclaim{Lemma~6.3} {\rm (\romannumeral1)} Let $n=2k+1${\rm,} then $f^\eta$
is nontrivial if and only if $\eta\ne0$.

{\rm (\romannumeral2)} Let $n=2k${\rm,} then $f^\eta$ is nontrivial if and 
only if $\eta\ne\bar\eta=(1,\dots,1,k\mod 2)$.
\endproclaim

Let now $\theta\:\V\to\K$ denote the projection on $\K$ along $\V_3$.
We lift $f^\eta$ to the whole $\V$ with the help of $\theta$ by defining
$F^\eta(X)=f^\eta(\theta(X))$ for any $X\in\V$.

It is easy to check that the standard coupling $(\cdot,\cdot)$ restricted to
$\Sigma_3$ identifies
$\K^*$ with $\V_3^\perp\subset\V^*$. Let $f^\eta=(f^\eta_1,\dots,f^\eta_k)
\in\V_3^\perp$; we define a linear functional 
$c^h=(c^h_1,\dots,c^h_k)\in\W^*$ as follows: $c^h=\Lambda f^{\psi_n(h)}$,
where $\Lambda$ is the matrix in $GL_k(\bF_2)$ such that $\lambda_{ij}=1$
iff $i\ls j$,
and lift it to the whole $T^n(\bF_2)$ with the help of $\tau$ by defining
$C^h(M)=c^h(\tau(M))$ for any $M\in T^n(\bF_2)$.

\proclaim{Lemma~6.4} Let $M\in T^n(\bF_2)$ be an arbitrary matrix in 
$S^h${\rm,} $h\in\bF_2^n$. If $\widetilde M\in S(M)$ and 
$\widetilde M\sim M${\rm,} then $C^h(\widetilde M)=C^h(M)$.
\endproclaim

\demo{Proof} Indeed, let $\eta=\psi_n(h)$ and let $M'=g_{ij}M$, $M'\ne M$. 
Then $C^h(M')=C^h(M)$ if 
$(i,j)\notin\Sigma_3$ and $C^h(M')-C^h(M)=c^h_j+c^h_{j+1}=f^{\eta}_j$
otherwise. Hence the total variation of $C^h(M)$ along the trajectory defined
by an arbitrary $M$-nonredundant word $\Omega$ equals 
$\var_\Omega C^h=[\Omega:\Supp F^{\eta}]\mod 2$. 
Therefore, to prove the Lemma it suffices to show that for any 
$M$-nonredundant word $\Omega$ such that $\widetilde M=\Omega M$ 
the number $[\Omega:\Supp F^{\eta}]$ is even.

Evidently, there exists a decomposition $\Omega=\Omega_1\Omega_2\dots\Omega_r$
with $\Omega_i\ne\varnothing$ for $i\gs 2$ such that all the elements of
$\Omega_{2i+1}$ belong to the subgroup $\G_n^1\subset\G_n$ generated by the 
elements of $\V_1$, and all the elements of $\Omega_{2i}$ belong to the 
subgroup $\G_n^2\subset\G_n$ generated by the elements of $\V_2$. By 
Lemma~4.5, instead of looking 
at the trajectory of $M_{2i}=\{\Omega_1\dots\Omega_{2i}\}\Psi_n(M)$ under
$\Omega_{2i+1}$ one can look at the trajectory of 
$(L^{\eta}_1)^{-1}(M_{2i})$ in $\V_1$, and instead of the trajectory of 
$M_{2i+1}=\Omega_{2i+1}M_{2i}$ under $\Omega_{2i+2}$, at the trajectory of 
$(L^{\eta}_2)^{-1}(M_{2i+1})$ in $\V_2$. These pieces are glued together 
with the help of the isomorphism $I^\eta$, see Fig.~7. Let us find the
total variation $\var_\Omega F^\eta$ of the functional $F^\eta$ along the 
obtained trajectory.

\vskip 10pt
\centerline{\hbox{\hskip 0.2cm \epsfysize=4cm\epsfbox{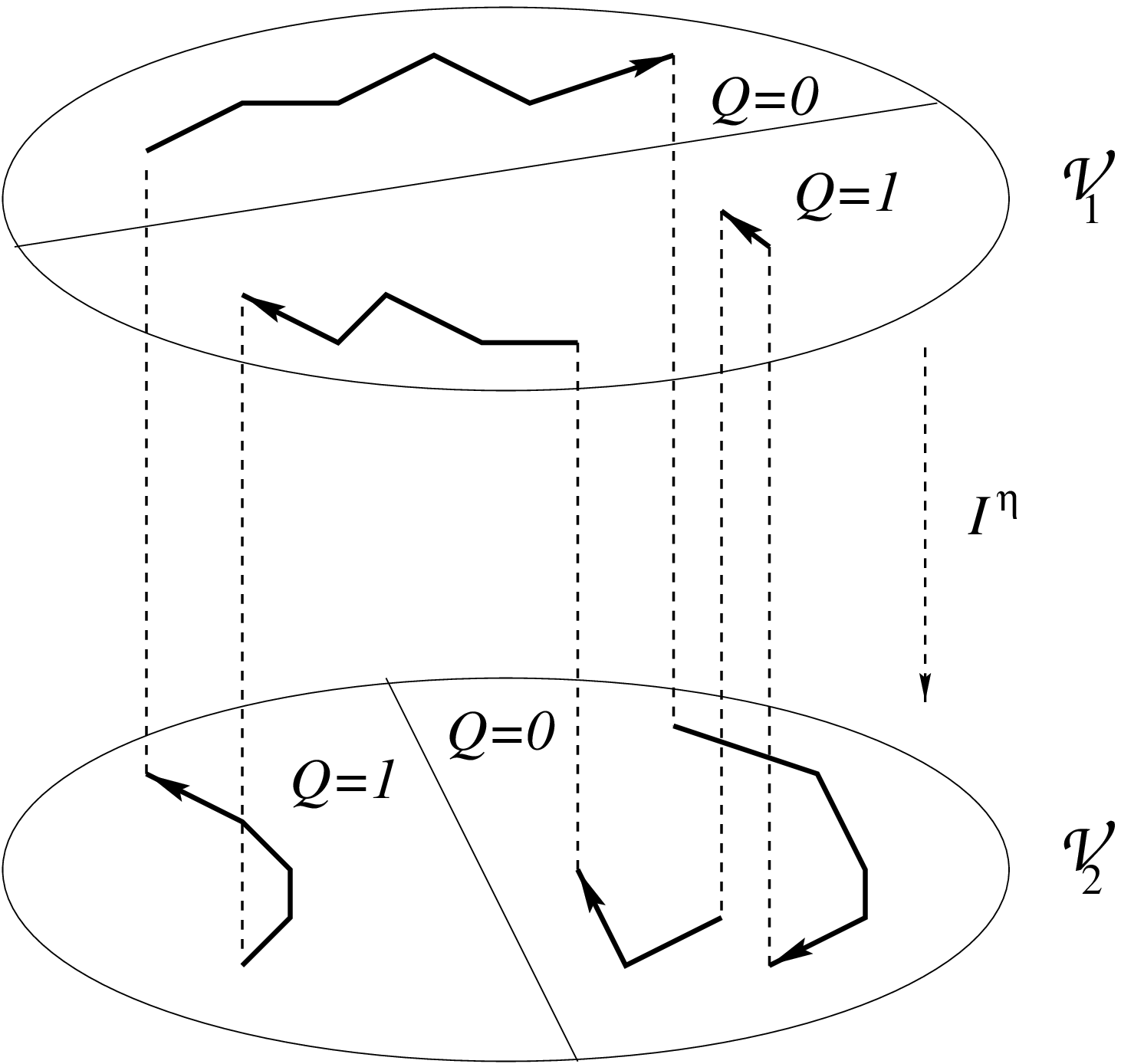}}}
\midspace{1mm} \caption{Fig.7. Trajectory defined by $\Omega$ in the
charts $\V_1$ and $\V_2$}
\vskip 5pt

It follows immediately from Lemma~3.2 that for any $X\in\V$ and $X'=g_{ij}X
\ne X$ one has $F^\eta(X')-F^\eta(X)=1$ iff $(i,j)\in\Supp F^\eta$. Therefore,
the total variation of $F^\eta$ along the pieces of the trajectory that
lie entirely in $\V_1$ or $\V_2$ equals $[\Omega:\Supp F^\eta]\mod 2$.

Let us find the variation of $F^\eta$ under the isomorphism $I^\eta$. 
Similarly to the proof of Lemma~5.9, for any $X\in\V_1$
one has $I^\eta(X)=X+H^\eta_3+X^\K$
with $X^\K\in\K$. Thus, $F^\eta(I^\eta(X))-F^\eta(X)=F^\eta(H^\eta_3)+
F^\eta(X^\K)=Q(X^\K)+(A^\eta_2,X^\K)$, since $\theta(H^\eta_3)=0$.

On the other hand, $Q(I^\eta(X))-Q(X)=Q(X^\K)+Q(H^\eta_3)+
(A^\eta_1+A^\eta_2,X)$. Besides, $(A^\eta_2,X)=0$ since $X\in\V_1$, 
$(A^\eta_2,X+H^\eta_3+X^\K)=0$ since $I^\eta(X))=X+H^\eta_3+X^\K\in\V_2$,
and $Q(H^\eta_3)+(A^\eta_2,H^\eta_3)=0$ by Lemma~6.2.
Adding the last four equalities we get $Q(I^\eta(X))-Q(X)=Q(X^\K)+
(A^\eta_2,X^\K)$. Therefore, $F^\eta(I^\eta(X))-F^\eta(X)=Q(I^\eta(X))-Q(X)$.

Recall that by Lemma~4.3 $Q$ is constant on each orbit in $\V_1$ and $\V_2$.
Therefore, $\var_\Omega F^\eta=\var_\Omega Q+[\Omega:\Supp F^\eta]\mod 2$.
Since $\Psi_n(M)=\Psi_n(\widetilde M)$, the trajectory defined by $\Omega$
is a loop, and hence $\var_\Omega F^\eta=\var_\Omega Q=0$, thus
$\var_\Omega C^h=[\Omega:\Supp F^{\eta}]\mod 2=0$.
\qed
\enddemo

Let us now prove the converse statement.

\proclaim{Lemma~6.5} Let $M\in T^n(\bF_2)\setminus\Cal I_n$ 
be an arbitrary matrix in $S^h${\rm,} $h\in\bF_2^n$. 
If $\widetilde M\in S(M)$ and $C^h(\widetilde M)=C^h(M)${\rm,} 
then $\widetilde M\sim M$.
\endproclaim

\demo{Proof} In what follows we assume that $h$ is fixed and $\eta=\psi_n(h)$.
First of all, for any $Z\in\K$ we choose two matrices $M_0^Z,M_1^Z\in\V_1$
such that
$$
\theta(I^\eta(M_i^Z))-\theta(M_i^Z)=Z,\qquad Q(M_i^Z)=i,\qquad i=0,1.
$$
Such a pair exists for any $Z\in\K$, since the first condition defines a 
translation of $\V_1\cap\V_2$, and $Q$ is nonlinear on $\V_1\cap\V_2$, and 
hence nonconstant on any translation of $\V_1\cap\V_2$. Moreover, we can 
further assume that $M_i^Z\ne0$ and $I^\eta(M_i^Z)\ne0$ for $i=0,1$.

Let us fix an arbitrary nontrivial matrix $X\in\V_1$. By Lemma~4.3, $X$,
$X'=M_{Q(X)}^Z$, and $X''=M_{Q(X)}^0$ belong to the same $\G_n^1$-orbit.
Let $\Omega'_Z$ be an $X$-nonredundant word in $\{g_{ij}\}$, 
$(i,j)\notin\Sigma_1$, such that $X'=\Omega'_ZX$. We denote $Y'=I^\eta(X')$,
$Y''=I^\eta(X'')$. By the proof of Lemma~6.4, the variation of $Q$ under 
$I^\eta$ equals the variation of $F^\eta$. Therefore, if $f^\eta(Z)=0$, then
$Q(Y')=Q(X')=Q(X)=Q(X'')=Q(Y'')$, and hence $Y',Y''\in\V_2$ belong to the
same $\G_n^2$-orbit. Let $\Omega''_Z$ denote a $Y'$-nonredundant word in
$\{g_{ij}\}$, $(i,j)\notin\Sigma_2$, such that $Y''=\Omega''_ZY'$. Evidently,
$\Omega_Z=\Omega'_Z\Omega''_Z\bar\Omega'_0$ is $X$-nonredundant and 
$\Omega_ZX=X$.

By Lemma~6.1, it suffices to prove that the value of any linear functional
other than $c^h$ can be changed along an appropriate $\Omega_Z$ with 
$f^\eta(Z)=0$. Take any $W\in\K$, and let
$f^W\in\V_3^\perp$ be the conjugate to $W$ with respect to $(\cdot,\cdot)$
restricted to $\Sigma_3$; put $F^W=f^W\circ\theta$. Then, similarly to the
proof of Lemma~6.4, $0=\var_{\Omega_Z}F^W=[\Omega_Z:\Supp F^W]+f^W(Z) \mod 2$.
Therefore, for $c^W=\Lambda f^W$ and $C^W=c^W\circ\tau$ one has 
$\var_{\Omega_Z}C^W= [\Omega_Z:\Supp F^W]\mod 2=f^W(Z)$. 

Let now $W=\sum_{i=1}^kw_i\widetilde P_{i}$. Evidently, it is 
enough to consider the following two types of $W$: (i) $w_i=0$ for $i\ne j$,
$w_j=1$, where $\eta_j=0$, and (ii) $w_i=0$ for $i\ne j_1,j_2$, 
$w_{j_1}=w_{j_2}=1$, where $\eta_{j_1}=\eta_{j_2}=1$. 
In order to get $\var_{\Omega_Z}C^W=1$ for $n=2k+1$ one takes
$Z=\widetilde P_{j}$ in the first case and $Z=\widetilde P_{j_1}+
\widetilde P_{j_3}$ with $\eta_{j_3}=1$ in the second case. For $n=2k$ one 
takes either $Z=\widetilde P_{j}+\widetilde P_{j_3}$ with $\eta_{j_3}=0$ 
or $Z=\widetilde P_{j_3}$ with $\eta_{j_3}=1$ in the first
case and $Z=\widetilde P_{j_1}+\widetilde P_{j_3}$ with $\eta_{j_3}=1$ in the 
second case.
\qed
\enddemo

\demo{Proof of Theorem~2.2} Let $n=2k+1$ and $\psi_n(h)\ne0$. 
Then by Theorem~2.5 and Lemmas~6.3,~6.4,~6.5 the stratum $S^h$ consists 
of two orbits distinguished by $C^h$. Their lengths are equal since $C^h$ 
is linear and nontrivial. For $\psi_n(h)=0$ Lemmas~6.3 and~6.5 imply that
nontrivial orbits in $S^h$ are just coverings of degree $2^k$ of nontrivial
orbits in $\Psi_n(S^h)=\D_{n-1}^\perp$ described in Theorem~2.5. Besides, 
each matrix in $\Cal I_n\cap S^h$ forms itself a trivial orbit. 
The case $n=2k$ is treated in the same way.
\qed\enddemo

\subheading{6.2. The structure of the symmetric strata of the first 
$\G_n$-action: an alternative approach} In \S 2.2 we have introduced 
the second $\G_n$-action as
the one induced by the first $\G_n$-action on the quotient of $T^n(\bF_2)$
modulo the subspace $\Cal I_n$ of the invariants of the first  $\G_n$-action.
Now we apply the same construction in the conjugate space. That is, we
consider the action induced by the conjugate to the first $\G_n$-action
on the  quotient of $(T^n(\bF_2))^*$ modulo the subspace $\DD_n$ of the 
dual invariants of the first  $\G_n$-action. One can suggest the following 
natural description of this induced action.

Consider the linear map $\Phi_n\: (T^n(\bF_2))^*\to  T^{n-1}(\bF_2)$ such
that the $(i,j)$th entry in the image equals the sum of the entries
$(i,j)$, $(i,j+1)$, $(i+1,j)$, and $(i+1,j+1)$ in the inverse image (as
before, entry $(i,j)$ in the image can be considered as representing the
submatrix $M^{ij}$ of the initial matrix). Lemma~3.1 implies immediately
that $\ker\Phi_n=\DD_n$, and we thus obtain the induced action of $\G_n$
on $T^{n-1}(\bF_2)$.

The following observation gives a clue to the structure of the symmetric 
strata of the first $\G_n$-action.

\proclaim{Lemma~6.6} The $\G_n$-action on $T^{n-1}(\bF_2)$ induced by $\Phi_n$
coincides with the second $\G_n$-action.
\endproclaim

\demo{Proof} Simple exercise in linear algebra. 
\qed
\enddemo

Therefore, the above construction provides an alternative description of the
second $\G_n$-action. As an immediate corollary of this description we obtain
the structure of the stratum $\DD_n^\perp$ at height $(0,\dots,0)$.

\proclaim{Lemma~6.7} The orbits of the first $\G_n$-action in the stratum
$\DD_n^\perp$ are isomorphic to the orbits of the conjugate to the second
$\G_n$-action.
\endproclaim

\demo{Proof} Indeed, consider the conjugate mapping $\Phi_n^*\: 
(T^{n-1}(\bF_2))^*\to T^n(\bF_2)$. It is easy to check that the image of
$\Phi_n^*$ coincides with $\DD_n^\perp$. Moreover, $\ker\Phi_n^*=0$, and
hence $\Phi_n^*$ provides an isomorphism between $(T^{n-1}(\bF_2))^*$ and
$\DD_n^\perp$. Now Lemma~6.6 implies that the first $\G_n$-action on
$\DD_n^\perp$ is isomorphic to the conjugate to the second $\G_n$-action.
\qed
\enddemo

 We thus get the part of Theorem~2.2 concerning the structure of symmetric 
strata.

\proclaim{Lemma~6.8} {\rm (\romannumeral1)} Let $n=2k+1\gs5${\rm,} then
each of $2^{k+1}$ symmetric strata consists of one orbit of length
$2^{2k^2+k-1}-\ee_k2^{k^2+k-1}${\rm,} one orbit of length 
$2^{2k^2+k-1}+\ee_k2^{k^2+k-1}-2^k${\rm,} and $2^k$ orbits of length $1${\rm,}
where $\ee_k=-1$ for $k=4t+1$ and $\ee_k=1$ otherwise{\rm;}

{\rm (\romannumeral2)} Let $n=2k\gs6${\rm,} then
each of $2^{k}$ symmetric strata consists of two orbits of length
$(2^{2k(k-1)}-1)2^{k-1}$ and $2^k$ orbits of length $1$.
\endproclaim

\demo{Proof}
 The additive group $\Cal I_n$ acts on $T^n(\bF_2)$ by 
translations. Evidently, any such translation takes $\DD_n^\perp$ to a 
symmetric stratum and any symmetric stratum is 
the translation of $\DD_n^\perp$ by a vector in $\Cal I_n$.
Besides, the action of $\Cal I_n$  commutes with
the first $\G_n$-action, and hence takes orbits to orbits. Therefore, the
orbit structure of any symmetric stratum coincides with that of the
stratum $\DD_n^\perp$.

To find the lengths of the orbits we use Lemmas~6.7,~5.6, and~4.3. 
Since $\ku=k$, we immediately get $2^k$ orbits of length $1$. 
The lengths of the other two orbits are obtained readily from 
Corollary~5.3 and Lemma~5.4. 
\qed
\enddemo

    \heading  \S  7. Final Remarks \endheading
\subheading{7.1} The technique of the previous paper \cite{SSV} is
substantially generalized in a forthcoming paper joint with A.~Zelevinsky to 
a wide class of intersections of pairs of Schubert cells and double Schubert 
cells.  This generalization is based  on the chamber
ansatz developed for the flag varieties and semisimple groups in 
\cite{BFZ, BZ, FZ}. For each reduced decomposition  of 
an element $u$ in a classical Coxeter group of type A, D, E we define a group 
which acts by symplectic transvections on the $\Bbb Z$-module 
generated by all chamber sets associated with the chosen decomposition. 
We prove that the $\bF_2$-reduction 
of the above action counts the number of connected components in 
in the intersection $B\cap B_u$ of two real open Schubert cells in $G/B$ 
taken in the split form. Analogous results are obtained for the 
intersections $G^{u,v}=BuB\cap B_{-}vB_{-}$ in a semisimple simplylaced 
group 
$G$.  
 
These results lead to the following general  setup. Given a connected 
undirected graph $\Gamma$, consider the $\bF_2$-vector space 
$\V_\Gamma$ generated by its vertices and define an $\bF_2$-valued bilinear 
form $\langle p,q\rangle=\sum\{p_iq_j\:\text{$i$ adjacent to $j$}\}$. 
Every vertex $\delta\in\Gamma$ determines the symplectic transvection 
$T_\delta:\V_\Gamma\to \V_\Gamma$ sending $p$ to $p-\langle p,\delta\rangle
\delta$. Let us also choose some subset $B$ of vertices of $\Gamma$ and 
define the group $\G_{B}$ generated by all $T_\delta$, $\delta\in B$.

\proclaim{Problem} Find the number of orbits of the $\G_{B}$-action on 
$\V_\Gamma$. 
\endproclaim

A substantial part of the results of the present paper are valid in this 
more general setup as well. The authors hope to solve the above Problem at 
least for the graphs (and their groups) arising from intersections of 
Schubert cells. Our preliminary considerations, supported by numerical
evidence, suggest the following surprising conjecture.

\proclaim{Conjecture} The number of connected components in the intersection
of two open Schubert cells in $SL_{n+1}(\Bbb R)/B$ in relative position $w$
equals $3\cdot 2^{n-1}$ for any generic $w$ and $n\gs5$.
\endproclaim

Here genericity means that the graph of bounded chambers introduced implicitly
in \cite{SSV} and studied in detail in the forthcoming paper with
A.~Zelevinsky contains an induced subgraph isomorphic to the Dynkin diagram
of $E_6$.
 
\subheading{7.2} Another very intriguing fact is that the group 
$\G_n$ studied in the present paper is the $\bF_2$-reduction of the 
monodromy group of the Verlinde algebra $su(3)_n$ introduced in 
\cite{Zu} and studied recently in \cite{GZV}.  (The authors are obliged 
to S.~Chmutov and S.~M.~Gusein-Zade for this observation.) Notice that 
 the monodromy group of any $su(m)_n$ is interpreted in \cite {GZV} as the 
monodromy group of the isolated singularity of a certain Newton polynomial 
expanded in elementary symmetric functions. This gives us a hope to both  
find a natural representation of Verlinde algebra $su(3)_n$ in the 
sections of some bundle over the intersection of open opposite 
Schubert cells and, more generally, to find a relation between the topology 
of intersections of Schubert cells and singularity theory of symmetric 
polynomials.

 \enddocument

\end
Observe that by Corollary~5.2(ii) we have $Q(K)=0$. Therefore, $Q(M)=0$ yields
$Q(I^\ga(M))=Q(H^\ga)+(M,a_1^\ga+a_2^\ga)$. Let $Q(H^\ga)=0$, then
to achieve $Q(I^\ga(M))=1$ one has to choose $M$ in such a way that
$(M,a_1^\ga+a_2^\ga)=1$. It follows that any $M$ such that $\Supp M$
consists of two nonadjacent vertices and intersects  
$\Supp a_1^\ga\cup\Supp a_2$ by exactly one vertex satisfies the conditions
$Q(M)=0$ and $Q(I^\ga(M))=1$ in this case.

Let now $Q(H^\ga)=1$, then  $Q(I^\ga(M))=1$ is equivalent to 
$(M,a_1^\ga+a_2^\ga)=0$. It follows that any $M$ such that $\Supp M$
consists of two nonadjacent vertices and does not intersect  
$\Supp a_1^\ga\cup\Supp a_2$ satisfies the conditions
$Q(M)=0$ and $Q(I^\ga(M))=1$ in this case.

Let us turn now to the choice of $N$. As above, $Q(N)=1$ yields
$Q(I^\ga(N))=1+Q(H^\ga)+(N,a_1^\ga+a_2^\ga)$. In the case $Q(H^\ga)=0$
any $N$ such that $\Supp N$ consists of one vertex not belonging to
$\Supp a_1^\ga\cup\Supp a_2$ satisfies the conditions
$Q(N)=1$ and $Q(I^\ga(N))=1$. In the case $Q(H^\ga)=1$ it suffices to take
any $N$ such that $\Supp N$ consists of one vertex belonging to
$\Supp a_1^\ga\cup\Supp a_2$. 

The existence of the matrices $M$ and $N$ as above in all the cases is evident,
since $\Supp a_i^\ga\subseteq S_i$, $i=1,2$, by the definition of $a_i^\ga$.
\qed
\enddemo

$$\matrix
&	\vdots	&
n=2k+1 & &	& \vdots
&	n=2k	&	&	& \vdots\\
\hdotsfor{10}\\
\text{type}& \vdots & \text{card}
& \sharp_{orb}&\text{mlt} &\vdots & \text{card} & \sharp_{orb}&\text{mlt}
&\vdots\\
\hdotsfor{10}\\
\text {trivial}& \vdots &	1	&1&2^{2k+1}&\vdots&1 &
1&2^{2k} &\vdots \\
\text {stand.}&\vdots &2^{2k^2}&2^k-1&2^{k+2}&\vdots &
2^{2k(k-1)}&2^{k}-2&2^{k+1}
&\vdots \\
\text {type 1}&\vdots &(2^{k^2}-1)2^{k^2-1}&1&2^{k+1}&\vdots & 0&0&0&\vdots \\
\text {type 2}&\vdots &(2^{k^2}-1)(2^{k^2-1}+1)&1&2^{k+1}&\vdots
&0&0&0&\vdots \\
\text {type 3}&\vdots &0&0&0&\vdots &
(2^{2k(k-1)}-1)&1&2^{k+1}&\vdots \\
\text {type 4}&\vdots &0&0&0&\vdots &
(2^{k(k-1)}-1)2^{k(k-1)-1}&1&2^{k}&\vdots \\ \text{type 5}&\vdots &0 &
0&0&\vdots &
(2^{k(k-1)}+1)2^{k(k-1)-1} &1&2^{k}&\vdots \\
\hdotsfor{10}
\endmatrix $$

\enddocument